\documentclass[times]{nlaauth}
\usepackage[mathscr]{euscript}
\usepackage{mathtools, paralist, stmaryrd}
\usepackage{array, multirow, cases, enumerate, float}
\usepackage{algorithm,algpseudocode}
\usepackage{hyphenat}
\usepackage{hhline}
\usepackage{xcolor,colortbl}
\usepackage[nice]{nicefrac}

\providecommand{\norm}[1]{\left\lVert#1\right\rVert}
\newcommand{\mtx}[1]{\ensuremath{\mathbf{#1}}}

\DeclareMathOperator{\tenrank}{rank}
\DeclareMathOperator{\numel}{numel}
\newcommand\Tstrut{\rule{0pt}{2.6ex}}

\graphicspath{{Images//}}

\usepackage{hyperref}

\begin{document}

\runningheads{H. De Sterck and A.J.M. Howse}{Nonlinear Acceleration by L\hyp{}BFGS}

\title{Nonlinearly Preconditioned L\hyp{}BFGS as an Acceleration Mechanism for Alternating Least Squares, with Application to Tensor Decomposition}

\author{Hans De Sterck\affil{1} and Alexander J.M. Howse\affil{2}\corrauth}

\address{\affilnum{1}School of Mathematical Sciences, Monash University, Melbourne, Australia
\break
\affilnum{2} Department of Applied Mathematics, University of Waterloo, Canada}

\corraddr{Department of Applied Mathematics, University of Waterloo, 200 University Ave W, Waterloo, ON, Canada, N2L 3G1. Email: {\tt ahowse@uwaterloo.ca}.}

\begin{abstract}
We derive nonlinear acceleration methods based on the limited memory BFGS (L\hyp{}BFGS) update formula for accelerating iterative optimization methods of alternating least squares (ALS) type applied to canonical polyadic (CP) and Tucker tensor decompositions. Our approach starts from linear preconditioning ideas that use linear transformations encoded by matrix multiplications, and extends these ideas to the case of genuinely nonlinear preconditioning, where the preconditioning operation involves fully nonlinear transformations. As such, the ALS-type iterations are used as fully nonlinear preconditioners for L\hyp{BFGS}, or, equivalently, L\hyp{BFGS} is used as a nonlinear accelerator for ALS. Numerical results show that the resulting methods perform much better than either stand-alone L\hyp{BFGS} or stand-alone ALS, offering substantial improvements in terms of time\hyp{}to\hyp{}solution and robustness over state\hyp{}of\hyp{}the\hyp{}art methods for large and noisy tensor problems, including previously described acceleration methods based on nonlinear conjugate gradients and nonlinear GMRES. Our approach provides a general L\hyp{}BFGS-based acceleration mechanism for nonlinear optimization.
\end{abstract}

\keywords{Nonlinear Acceleration; Nonlinear Optimization; Nonlinear Preconditioning; Tensor Decompositions; Canonical Polyadic Decomposition; Tucker Decomposition; Quasi-Newton Methods; L\hyp{}BFGS}

\maketitle

\section{Introduction}

Simple nonlinear iterative optimization methods like alternating least squares (ALS) and coordinate descent (CD) are widely used in a variety of application domains, including tensor decomposition, image processing, computational statistics and machine learning \cite{kolda2009tensor,wright2015coordinate}.
These methods solve large-scale unconstrained nonlinear optimization problems 
\begin{equation}
\min \quad f(\mtx{x}),
\label{eq:opt}
\end{equation}
by fixing, in each iteration, most components of the variable vector x at their values from the current iteration, and approximately or exactly minimizing $f(\mtx{x})$ with respect to the remaining components.
These simple methods are competitive with more sophisticated approaches in a variety of contexts \cite{kolda2009tensor,wright2015coordinate}.

As an example in point, methods of ALS type are workhorse algorithms for canonical polyadic (CP) or Tucker tensor decompositions \cite{kolda2009tensor,acar2011scalable}, which have applications in multimodal data analysis and compression. 
For the Tucker decomposition, the standard alternating algorithm is the higher-order orthogonal iteration (HOOI) \cite{de2000best}.
It is not an ALS method per se, but is of a similar nature: it maximizes in an alternating fashion quadratic objectives that are obtained by fixing most components of the variable to be optimized. For ease of exposition we will refer to HOOI as an ALS-type method.
These methods of ALS type, for CP and Tucker, are widely used, but by themselves may converge too slowly when problems become ill-conditioned or when accurate solutions are required, even though they are under those circumstances often still more efficient than more advanced alternatives like the limited memory Broyden\hyp{}Fletcher\hyp{}Goldfarb\hyp{}Shanno (L\hyp{}BFGS) quasi\hyp{}Newton method or nonlinear conjugate gradient (NCG) method \cite{acar2011scalable}.

For these situations where convergence is slow, nonlinear acceleration mechanisms for the simple optimization methods are of significant interest.

In previous work it was shown that the NCG method and a nonlinear GMRES (NGMRES) optimization method can be used as nonlinear convergence accelerators for ALS-type methods in tensor decomposition problems, substantially improving convergence and leading to state-of-the-art methods for tensor decomposition in terms of performance \cite{sterck2012nonlinear, sterck2014nonlinearly, sterck2016nonlinearly}.

However, as the most efficient quasi-Newton (QN) method for optimization, the performance of L\hyp{}BFGS is expected to be superior to NCG and NGMRES, and this naturally suggests to use L\hyp{}BFGS as a nonlinear convergence accelerator, with the promise of improved performance compared to acceleration by NCG or NGMRES. In this paper we derive two approaches for using L\hyp{}BFGS as a nonlinear convergence accelerator for ALS-type optimization methods, and show numerically that nonlinear acceleration of ALS-type methods for tensor decomposition by L\hyp{}BFGS indeed leads to substantial performance gains compared to other leading methods for tensor decomposition, when problems are ill-conditioned or accurate solutions are required.

The derivation of our method for nonlinear acceleration by L\hyp{}BFGS is rooted in the case of minimizing convex quadratic functionals that correspond to solving symmetric positive definite (SPD) linear systems. It has long been known that iteration growth resulting from ill-conditioning can in this case be combated in highly efficient manners by preconditioning the problem, where linear transformations are used to improve conditioning and reduce iteration counts. Preconditioning using linear transformations encoded by matrix multiplications has been generalized to NCG and L\hyp{}BFGS beyond the convex quadratic case, and the resulting linearly preconditioned methods for nonlinear optimization are well-known in the optimization community \cite{hager2006survey, luenberger2015linear}. However, in order to use L\hyp{}BFGS as a nonlinear acceleration mechanism for methods such as ALS, one has to extend the notion of preconditioning by linear transformations to the case of genuinely nonlinear preconditioning, where the preconditioning operation involves fully nonlinear transformations, as opposed to linear transformations encoded by matrix multiplications. This was done in \cite{sterck2012nonlinear, sterck2014nonlinearly, sterck2016nonlinearly} for the NCG and NGMRES methods, in ways that reduce to the linearly preconditioned case for convex quadratic objectives, and it is our goal to do this in this paper for L\hyp{}BFGS. When applying the resulting methods, e.g., to ALS-type iterations for tensor decomposition, one can view ALS as the fully nonlinear preconditioner for L\hyp{}BFGS, or, equivalently, one can view L\hyp{}BFGS as a nonlinear acceleration mechanism for ALS.

These methods derived from nonlinear preconditioning ideas can be situated in the broader context of the current extensive research interest in nonlinear preconditioning for solving systems of nonlinear equations $\mtx{g}(\mtx{x}) = \mtx{0}$ \cite{cai2002nonlinearly}, as summarized in the review paper by Brune et al. \cite{brune2015composing}. Ideas on nonlinear preconditioning date back to the 1960s \cite{anderson1965iterative, bartels1974conjugate, concus1978numerical}, but they remain underexplored in theory and in practice \cite{brune2015composing}, especially in comparison to linear preconditioning. We extend some of the strategies used in \cite{brune2015composing} and in other papers on nonlinear preconditioning \cite{sterck2012nonlinear, sterck2014nonlinearly, sterck2016nonlinearly} to L\hyp{}BFGS in the optimization context. We also derive nonlinearly preconditioned versions of the limited memory Broyden (L\hyp{}Broyden) QN method. Without preconditioning, L\hyp{}Broyden is often inferior to L\hyp{}BFGS, but, interestingly, we find that, with nonlinear preconditioning, L\hyp{}Broyden can sometimes be competitive with leading methods.

We consider two strategies for extending nonlinear preconditioning to L\hyp{}BFGS. The first approach is based on the general idea of left preconditioning, which is commonly used when solving linear systems, and has been generalized to the case of nonlinear preconditioners for nonlinear systems of equations, see, e.g., the review paper \cite{brune2015composing}. We refer to this approach as the LP approach (for Left Preconditioning). Nonlinear left preconditioning of L\hyp{}BFGS has been considered before in \cite{brune2015composing} in the context of nonlinear systems solvers for partial differential equations (PDEs), and in this paper we investigate how the idea of left preconditioning applies to the optimization context, using ALS-type preconditioners.

Our second approach is inspired by a linear change of variables $\mtx{x}=\mtx{Cz}$ in optimization problem (\ref{eq:opt}). It is derived by applying the optimization method to $\widehat{f}(\mtx{z})=f(\mtx{Cz})$, and transforming back to the $\mtx{x}$ variables. The resulting linearly preconditioned methods are well-known in the optimization community \cite{hager2006survey, luenberger2015linear}, but in this paper we generalize this approach to the case of fully nonlinear preconditioners using nonlinear transformations (i.e., not encoded by matrix multiplications). We refer to this second approach as the TP approach (for Transformation Preconditioning). The resulting TP form of nonlinearly preconditioned L\hyp{}BFGS has, to our knowledge, not been considered before.

Both the LP and TP approaches were explored implicitly in the previous work on nonlinearly preconditioned NCG optimization methods \cite{sterck2014nonlinearly, sterck2016nonlinearly}. In that context, the two approaches arise naturally and have a simple interpretation. However, in the case of L\hyp{}BFGS, the situation is much more intricate, and in this paper we present a framework that explicitly formalizes the LP and TP approaches for general optimization methods, which allows us to derive nonlinearly preconditioned versions of the L\hyp{}BFGS and L\hyp{}Broyden methods.

We first illustrate the main ideas in the simplified context of linear preconditioning for convex quadratic optimization problems that correspond to solving SPD linear systems. This allows us to explain and interpret the LP and TP approaches in relation to well-known ideas for linear preconditioning. Some preliminary numerical tests will illustrate the two preconditioning approaches and will provide some initial numerical justification. We then extend the formalism to nonlinear preconditioning for general nonlinear objective functions $f(\mtx{x})$, and provide extensive numerical tests illustrating and comparing the merits of the two approaches.

To demonstrate the efficacy of the proposed nonlinearly preconditioned quasi\hyp{}Newton (NPQN) algorithms we use them to solve approximate tensor decomposition problems: decomposition of a multidimensional array into a sum and/or product of multiple components (for example, vector outer products) to reduce storage costs and simplify further data analysis. Such problems are typically cast as minimizing the approximation error in a given norm. Specifically, we consider decompositions into CP and Tucker tensor formats, the former representing a tensor as a sum of rank\hyp{}one terms and the latter as a multilinear tensor\hyp{}matrix product. Background material on computing tensor decompositions are presented in Appendix~\ref{appx:tensors}. Both of these tensor formats have standard iterative algorithms of ALS type for computing decompositions, which can be slow to converge when used independently, but can be useful as nonlinear preconditioners, resulting in significant acceleration. The resulting methods combine L\hyp{}BFGS with ALS, and are much more efficient than either of ALS or L\hyp{}BFGS by themselves, offering substantial improvements in terms of time\hyp{}to\hyp{}solution and robustness over state\hyp{}of\hyp{}the\hyp{}art methods for large and noisy tensor problems. 

As stated before, NGMRES \cite{sterck2012nonlinear} and NCG \cite{sterck2014nonlinearly} methods were used before to accelerate ALS for the CP decomposition problem, where the resulting combined methods are currently among the fastest available for noisy problems and when high accuracy is required.
The nonlinearly preconditioned NCG (NPNCG) and nonlinearly preconditioned NGMRES (NPNGMRES) algorithms were extended to the Tucker decomposition problem in \cite{sterck2016nonlinearly}, which involved the use of matrix manifold optimization to handle lack of uniqueness and equality constraints required for decomposition into a particular subtype of Tucker tensor called the higher order singular value decomposition (HOSVD). The details of adapting our proposed NPQN algorithms to the required manifold structure are provided in Appendix~\ref{appx:manifolds}. Several other matrix manifold optimization strategies for Tucker decompositions have been proposed. Newton's method was considered by Eld{\'e}n and Savas in \cite{elden2009newton}, followed by adaptation of the BFGS QN method, as well as L\hyp{}BFGS, by Savas and Lim \cite{savas2010quasi}. Manifold NCG \cite{ishteva2009numerical}, a Riemannian trust-region method \cite{ishteva2011best}, and a differential-geometric Newton's method \cite{ishteva2009differential} have been developed by Ishteva et al. In \cite{sterck2016nonlinearly}, the NGMRES- and NCG-accelerated ALS-type approaches were compared with these methods, and showed substantial improvements in performance for difficult problems.

The remainder of the paper is organized as follows. In Section~\ref{sect:optimization:methods} we describe the optimization methods that are used as building blocks in subsequent discussion. Section~\ref{sect:lin:precond} recalls linear preconditioning strategies in the context of minimizing quadratic objective functions which guide our nonlinear method development, and in Section~\ref{sect:nonlinear:preconditioning} we describe our new QN nonlinear preconditioning/nonlinear acceleration strategies, including how such methods can be informed by the linear case. In Section~\ref{sect:numerical:results} we provide numerical results for tensor decomposition problems, illustrating the improvements offered by the proposed methods, and finally in Section~\ref{sect:conclusions} we summarize our results and discuss future work.

\section{Optimization Method Building Blocks}
\label{sect:optimization:methods}

\subsection{Nonlinear Conjugate Gradient Method}
\label{subsect:conjugate:gradient}

The conjugate gradient (CG) method \cite{nocedal2006numerical} is an iterative solver for linear systems $\mtx{Ax} = \mtx{b}$ with symmetric positive definite (SPD) matrices $\mtx{A}$; or equivalently, a solver that minimizes convex quadratic objective functions 
\begin{equation}
f(\mtx{x})=\tfrac{1}{2}\mtx{x}^\intercal\mtx{A}\mtx{x}-\mtx{b}^\intercal\mtx{x}.
\label{eq:quadratic:objective}
\end{equation}
Starting from initial guess $\mtx{x}_0$ with initial search direction and residual $\mtx{p}_0 = \mtx{r}_0 = \mtx{b} - \mtx{A}\mtx{x}_0$. CG is based on the conjugacy of the sequence of search directions $\mtx{p}_k$ with respect to $\mtx{A}$, and generates an orthogonal sequence of residual vectors $\mtx{r}_k$ \cite{nocedal2006numerical}. In addition to the low storage requirements, we only require the means to compute the matrix\hyp{}vector product $\mtx{A}\mtx{p}_k$; storage of $\mtx{A}$ is unnecessary.


The nonlinear conjugate gradient (NCG) iteration (Algorithm \ref{alg:NCG}) arose as an adaptation of CG for minimizing general nonlinear objective functions $f(\mtx{x})$ \cite{nocedal2006numerical}. NCG generalizes CG as follows: 
\begin{inparaenum}[(i)]
\item we replace the residual $\mtx{r}_k$ with the gradient of the objective function, $\mtx{g}(\mtx{x}) = \nabla f(\mtx{x})$,
\item the step length $\alpha_k$ must be determined by a line\hyp{}search, and
\item the search direction update parameter $\beta_k$ can be specified by a number of different formulas.
\end{inparaenum}
Three of the most successful are
\begin{align}
\label{eq:beta:PR}
\text{Polak-Ribi\`ere \cite{polak1969note}:} \quad \beta_{k}^{\text{PR}} & = \frac{\mtx{g}_{k}^\intercal\mtx{y}_{k-1}}{\mtx{g}_{k-1}^\intercal\mtx{g}_{k-1}}, \\
\label{eq:beta:HS}
\text{Hestenes-Stiefel \cite{hestenes1952methods}:} \quad \beta_{k}^{\text{HS}} & = \frac{\mtx{g}_{k}^\intercal\mtx{y}_{k-1}}{\mtx{p}_{k-1}^\intercal\mtx{y}_{k-1}}, \\
\text{Hager-Zhang \cite{hager2005new}:} \quad \label{eq:beta:HZ}
\beta_{k}^{\text{HZ}} & = \left(\mtx{y}_{k-1} - 2\mtx{p}_{k-1}\frac{\norm{\mtx{y}_{k-1}}^2}{\mtx{p}_{k-1}^\intercal\mtx{y}_{k-1}}\right)^\intercal \frac{\mtx{g}_{k}}{\mtx{p}_{k-1}^\intercal\mtx{y}_{k-1}},
\end{align}
where $\mtx{g}_k = \mtx{g}(\mtx{x}_k)$ and $\mtx{y}_{k-1} = \mtx{g}_{k} - \mtx{g}_{k-1}$. One iteration of NCG is described in Algorithm \ref{alg:NCG}. Like the linear version, NCG enjoys the benefits of very low storage requirements.

\begin{algorithm}[ht]
\caption{Nonlinear Conjugate Gradients Iteration}\label{alg:NCG}
\begin{algorithmic}[1]
\Procedure{NCG}{$\mtx{g}(\cdot),\mtx{x}_k,\mtx{p}_{k},\mtx{g}_{k}$}
\State $\mtx{x}_{k+1} = \mtx{x}_k + \alpha_k\mtx{p}_k$ \Comment $\alpha_k$ determined by line\hyp{}search
\State $\mtx{g}_{k+1} = \mtx{g}(\mtx{x}_{k+1})$
\State Compute $\beta_{k}$ by one of (\ref{eq:beta:PR}--\ref{eq:beta:HZ})
\State $\mtx{p}_{k+1} = -\mtx{g}_{k+1} + \beta_{k}\mtx{p}_{k}$
\State\Return $\mtx{x}_{k+1},\,\mtx{p}_{k+1},\,\mtx{g}_{k+1}$
\EndProcedure
\end{algorithmic}
\end{algorithm}

In the linear case, the $\alpha_k$ and $\beta_k$ parameters are given by
$$
\alpha_k = (\mtx{r}_k^\intercal\mtx{r}_k)/(\mtx{p}_k^\intercal\mtx{A}\mtx{p}_k)
$$
and
$$
\beta_{k} = (\mtx{r}_{k+1}^\intercal\mtx{r}_{k+1})/(\mtx{r}_{k}^\intercal\mtx{r}_{k}).
$$

\subsection{Limited Memory QN Methods}
\label{subsect:limited:memory:QN}

QN methods \cite{dennis1996NMU} are iterations based on the standard Newton\hyp{}Raphson iteration for solving nonlinear systems $\mtx{g}(\mtx{x})$ or minimizing nonlinear functions $f(\mtx{x})$, depending on the context. The expensive evaluation of the Jacobian or Hessian matrix at each iteration is replaced with a low\hyp{}rank update of a matrix approximation based on a secant condition. The quadratic rate of convergence of Newton's method is traded for super\hyp{}linear convergence, with the hope that the approximation results in a significantly lower per iteration cost.

Limited memory QN iterations aim for further savings in storage requirements and work per iteration by expressing the matrix approximation in terms of an initial matrix (often diagonal) and at most $m$ vector pairs. In this way the full matrix approximation does not need to be formed or stored to compute matrix\hyp{}vector products. We describe two limited memory QN methods in the remainder of this subsection: one based on the good Broyden update for solving nonlinear systems (L\hyp{}Broyden), and one based on the BFGS update for minimizing nonlinear objective functions (L\hyp{}BFGS).

\subsubsection{The L\hyp{}Broyden Update}

We first describe the general good Broyden update for standard QN methods, and then give the limited memory variant. When considering a nonlinear system $\mtx{g}(\mtx{x}) = \mtx{0}$, we denote an approximation to the Jacobian matrix by $\mtx{A}_k$, and define the vectors $\mtx{s}_k = \mtx{x}_{k+1} - \mtx{x}_{k}$ and $\mtx{y}_k = \mtx{g}_{k+1} - \mtx{g}_k$ \cite{dennis1996NMU}. Broyden's good update minimizes the change in the affine model
$$
\mtx{M}_{k+1}(\mtx{x}) = \mtx{g}_{k+1} + \mtx{A}_{k+1}(\mtx{x}-\mtx{x}_{k+1})
$$
between iterations, subject to the secant equation
$$
\mtx{A}_{k+1}\mtx{s}_k = \mtx{y}_k.
$$
The resulting rank\hyp{}one update is
\begin{equation}
\mtx{A}_{k+1} = \mtx{A}_k + \dfrac{(\mtx{y}_k - \mtx{A}_k\mtx{s}_k)\mtx{s}_k^\intercal}{\mtx{s}_k^\intercal \mtx{s}_k},
\label{eq:broyden:update}
\end{equation}
which, by applying the Sherman-Morrison-Woodbury formula \cite{sherman1950adjustment, woodbury1950inverting, bartlett1951inverse}, gives the inverse matrix update
\begin{equation}
\mtx{A}_{k+1}^{-1} = \mtx{A}_k^{-1} + \dfrac{(\mtx{s}_k - \mtx{A}_k^{-1}\mtx{y}_k)\mtx{s}_k^\intercal\mtx{A}_k^{-1}}{\mtx{s}_k^\intercal\mtx{A}_k^{-1}\mtx{y}_k}.
\label{eq:inverse:broyden:update}
\end{equation}
An example of one QN iteration in this context is given in Algorithm \ref{alg:QN}.

\begin{algorithm}[ht]
\caption{Quasi-Newton Iteration for Nonlinear Systems}\label{alg:QN}
\begin{algorithmic}[1]
\Procedure{QN}{$\mtx{g}, \mtx{x}_k$, $\mtx{A}_{k}^{-1}$}
\State $\mtx{p}_k = -\mtx{A}_k^{-1}\mtx{g}(\mtx{x}_k)$
\State $\mtx{x}_{k+1} = \mtx{x}_k + \alpha_k\mtx{p}_k$ \Comment $\alpha_k$ determined by line\hyp{}search
\State $\mtx{s}_k = \mtx{x}_{k+1} - \mtx{x}_k$
\State $\mtx{y}_k = \mtx{g}(\mtx{x}_{k+1}) - \mtx{g}(\mtx{x}_k)$
\State $\mtx{A}_{k+1}^{-1} = \mtx{U}(\mtx{A}_{k}^{-1}, \mtx{s}_{k}, \mtx{y}_{k})$ \Comment $\mtx{U}$ update formula \eqref{eq:inverse:broyden:update}.
\State\Return $\mtx{x}_{k+1},\,\mtx{A}_{k+1}^{-1}$
\EndProcedure
\end{algorithmic}
\end{algorithm}

In the limited memory context, where only a window of $m$ previous vector pairs are retained, a compact, non\hyp{}recursive representation of update \eqref{eq:inverse:broyden:update} derived in \cite{byrd1994representations} is
\begin{equation}
\mtx{A}_k^{-1} = \left[\mtx{A}_0^{(k)}\right]^{-1} - \left(\left[\mtx{A}_0^{(k)}\right]^{-1}\mtx{Y}_k - \mtx{S}_k\right)\left(\mtx{M}_k + \mtx{S}_k^\intercal\left[\mtx{A}_0^{(k)}\right]^{-1}\mtx{Y}_k\right)^{-1}\mtx{S}_k^\intercal\left[\mtx{A}_0^{(k)}\right]^{-1}
\label{eq:inverse:broyden:compact}
\end{equation}
where 
\begin{align}
\mtx{S}_k & = [\mtx{s}_{k-m}\,|\, \mtx{s}_{k-m+1} \,|\, \cdots \,|\, \mtx{s}_{k-1}], \label{eq:defn:Sk}\\
\mtx{Y}_k & = [\mtx{y}_{k-m}\,|\, \mtx{y}_{k-m+1} \,|\, \cdots \,|\, \mtx{y}_{k-1}], \label{eq:defn:Yk}\\
\intertext{and}
(\mtx{M}_k)_{i,j} & = 
\begin{cases}
-\mtx{s}_{i-1}^\intercal\mtx{s}_{j-1} & \text{ if } i > j \\
0 & \text{ otherwise}
\end{cases}.\label{eq:defn:Mk}
\end{align}
For the initial Jacobian approximation $\mtx{A}_0^{(k)}$ we will typically use a scaled identity matrix, $\mtx{A}_0^{(k)}=\eta_k \mtx{I}$. It is this representation that we will use in the L\hyp{}Broyden iteration.

\subsubsection{The L\hyp{}BFGS Update}

In the optimization context we use $\mtx{B}_k$ to denote an approximation to the Hessian and $\mtx{H}_k$ for an approximation to the inverse of the Hessian. Arguably the most successful QN update for nonlinear optimization is the BFGS update, which, in addition to enforcing the secant equation, ensures that $\mtx{B}_k$ is SPD provided $\mtx{s}_{k-1}^\intercal\mtx{y}_{k-1} > 0$ and $\mtx{B}_{k-1}$ is SPD \cite{dennis1996NMU}. This is a rank\hyp{}two update given by
\begin{equation}
\mtx{B}_{k+1} = \mtx{B}_{k} + \dfrac{\mtx{y}_k\mtx{y}_k^\intercal}{\mtx{y}_k^\intercal\mtx{s}_k} - \dfrac{\mtx{B}_k\mtx{s}_k(\mtx{B}_k\mtx{s}_k)^\intercal}{\mtx{s}_k^\intercal\mtx{B}_k\mtx{s}_k},
\label{eq:bfgs:update}
\end{equation}
with inverse update
\begin{equation}
\mtx{H}_{k+1} = \mtx{H}_k + \dfrac{(\mtx{s}_k - \mtx{H}_k\mtx{y}_k)\mtx{s}_k^\intercal + \mtx{s}_k(\mtx{s}_k - \mtx{H}_k\mtx{y}_k)^\intercal }{\mtx{y}_k^\intercal\mtx{s}_k} - \dfrac{\langle\mtx{s}_k - \mtx{H}_k\mtx{y}_k, \mtx{y}_k\rangle \mtx{s}_k\mtx{s}_k^\intercal }{(\mtx{y}_k^\intercal\mtx{s}_k)^2}.
\label{eq:inverse:bfgs:update}
\end{equation}

In the limited memory case, given initial Hessian approximation $\mtx{B}_0^{(k)}$ $\left(\text{or }\mtx{H}_0^{(k)}\right)$ and at most $m$ vector pairs $\mtx{y}_k$, $\mtx{s}_k$, compact versions of \eqref{eq:bfgs:update} and \eqref{eq:inverse:bfgs:update} are \cite{byrd1994representations}:
\begin{equation}
\mtx{B}_k = \mtx{B}_0^{(k)} - 
\begin{bmatrix}
\mtx{B}_0^{(k)}\mtx{S}_k & \mtx{Y}_k
\end{bmatrix}
\begin{bmatrix}
\mtx{S}_k^\intercal \mtx{B}_0^{(k)} \mtx{S}_k & \mtx{L}_k \\
\mtx{L}_k^\intercal & -\mtx{D}_k
\end{bmatrix}
\begin{bmatrix}
\mtx{S}_k^\intercal \mtx{B}_0^{(k)} \\
\mtx{Y}_k^\intercal
\end{bmatrix}
\label{eq:bfgs:compact}
\end{equation}
and
\begin{equation}
\mtx{H}_k = \mtx{H}_0^{(k)} + 
\begin{bmatrix}
\mtx{S}_k & \mtx{H}_0^{(k)}\mtx{Y}_k
\end{bmatrix}
\begin{bmatrix}
\mtx{R}_k^{-\intercal}(\mtx{D}_k \mtx{Y}_k^\intercal \mtx{H}_0^{(k)} \mtx{Y}_k)\mtx{R}_k^{-1} & -\mtx{R}_k^{-\intercal}  \\
-\mtx{R}_k^{-1}  & \mtx{0}
\end{bmatrix}
\begin{bmatrix}
\mtx{S}_k^\intercal \\
\mtx{Y}_k^\intercal \mtx{H}_0^{(k)}
\end{bmatrix},
\label{eq:inverse:bfgs:compact}
\end{equation}
where $\mtx{S}_k$ and $\mtx{Y}_k$ are as in \eqref{eq:defn:Sk} and \eqref{eq:defn:Yk}, and
\begin{align}
\mtx{D}_k & = \text{diag}[\mtx{s}_{k-m}^\intercal\mtx{y}_{k-m}, \ldots, \mtx{s}_{k-1}^\intercal\mtx{y}_{k-1}], \label{eq:defn:Dk}\\
(\mtx{L}_k)_{i,j} & = \begin{cases}
(\mtx{s}_{k-m-1+i})^\intercal(\mtx{y}_{k-m-1+j}) & \text{ if } i > j \\
0 & \text{ otherwise }
\end{cases}, \label{eq:defn:Lk}\\
(\mtx{R}_k)_{i,j} & = \begin{cases}
(\mtx{s}_{k-m-1+i})^\intercal(\mtx{y}_{k-m-1+j}) & \text{ if } i \leq j \\
0 & \text{ otherwise }
\end{cases}. \label{eq:defn:Rk}
\end{align}
It is common to set $\mtx{H}_0^{(k)} = \gamma_k\mtx{I}$, where
\begin{equation}
\gamma_k = \frac{\mtx{s}_{k-1}^\intercal\mtx{y}_{k-1}}{\mtx{y}_{k-1}^\intercal\mtx{y}_{k-1}}.
\label{eq:defn:gamma}
\end{equation}
This choice of $\gamma_k$ is a scaling factor which attempts to make the size of $\mtx{H}_0^{(k)}$ similar to that of the true Hessian inverse $\nabla^2 f(\mtx{x}_{k-1})^{-1}$ along the most recent search direction, which helps ensure the search direction $\mtx{p}_k$ is scaled so that a unit step length $\alpha_k$ is acceptable in more iterations \cite{nocedal2006numerical}. 
(This is, in fact, a rudimentary form of preconditioning, using a simple linear transformation, on which we will improve below.)
When working with the inverse L\hyp{}BFGS update, the product $-\mtx{H}_k\mtx{g}_k$ defining the QN direction $\mtx{p}_k$ (similar to Algorithm \ref{alg:QN}) can be efficiently computed by a two\hyp{}loop recursion, described in Algorithm \ref{alg:2loop}.
The L\hyp{}BFGS update formula is then given by
$$
\mtx{x}_{k+1} = \mtx{x}_k + \alpha_k\mtx{p}_k,
$$
where $\alpha_k$ is determined by a line\hyp{}search.

\begin{algorithm}[ht]
\caption{L\hyp{}BFGS Two\hyp{}Loop Recursion}\label{alg:2loop}
\begin{algorithmic}[1]
\Procedure{2LOOP}{$\mtx{H}_0^{(k)}, \mtx{g}_k, \mtx{S}_k, \mtx{Y}_k$}
\State $\mtx{q} = \mtx{g_k}$
\For{$i=k-1,k-2,\ldots,k-m$}
\State $\rho_i =(\mtx{y}_i^\intercal\mtx{s}_i)^{-1}$
\State $\alpha_i = \rho_i \mtx{s}_i^\intercal\mtx{q}$
\State $\mtx{q} = \mtx{q} -\alpha_i\mtx{y}_i$
\EndFor
\State $\mtx{r} = \mtx{H}_0^{(k)}\mtx{q}$
\For{$i=k-m,k-m+1,\ldots,k-1$}
\State $\beta = \rho_i \mtx{y}_i^\intercal\mtx{r}$
\State $\mtx{r} = \mtx{r} + (\alpha_i-\beta)\mtx{s}_i$
\EndFor
\State\Return $\mtx{r}$ \Comment{Contains $\mtx{H_k}\mtx{g}_k$}
\EndProcedure
\end{algorithmic}
\end{algorithm}

In situations where $\mtx{s}_{k-1}^\intercal\mtx{y}_{k-1} \leq 0$, there is a damped BFGS variant that ensures the updated Hessian is SPD \cite{nocedal2006numerical} by defining
$$
\theta_k = \begin{cases}
1 & \text{ if } \mtx{s}_{k}^\intercal\mtx{y}_{k} \geq 0.1\mtx{s}_{k}^\intercal\mtx{B}_{k}\mtx{s}_{k} \\
(0.9\mtx{s}_{k}^\intercal\mtx{B}_{k}\mtx{s}_{k})/(\mtx{s}_{k}^\intercal\mtx{B}_{k}\mtx{s}_{k} - \mtx{s}_k^\intercal\mtx{y}_k) & \text{ if } \mtx{s}_{k}^\intercal\mtx{y}_{k} < 0.1\mtx{s}_{k}^\intercal\mtx{B}_{k}\mtx{s}_{k}
\end{cases}
$$
and setting
$$
\mtx{y}_k = \theta_k \mtx{y}_k + (1-\theta_k)\mtx{B}_k\mtx{s}_k,
$$
which reduces to the standard update for $\mtx{s}_{k}^\intercal\mtx{y}_{k} \geq 0.1\mtx{s}_{k}^\intercal\mtx{B}_{k}\mtx{s}_{k}$. We use this damping step in our L\hyp{}BFGS implementations.

\subsubsection{Relationship of BFGS to CG}
\label{subsubsect:bfgs:cg}

There are some noteworthy similarities between the CG and BFGS methods being considered, both for the convex quadratic objective function \eqref{eq:quadratic:objective} and more general nonlinear objective functions. It has been shown for convex quadratic objective functions that the CG and BFGS iterations are identical when exact line\hyp{}searches are used \cite{buckley1978extending, nazareth1979relationship}. Furthermore, the ``memory\hyp{}less'' BFGS method (L\hyp{}BFGS with $m=1$), in conjunction with an exact line\hyp{}search, applied to a general nonlinear objective function is equivalent to using NCG with the Hestenes\hyp{}Stiefel (HS) or Polak\hyp{}Ribi\`ere (PR) $\beta$ formulas (which are equivalent since $\mtx{g}_{k+1}^\intercal\mtx{p}_k = 0$ by the exact line\hyp{}search) \cite{nocedal2006numerical, luenberger2015linear}.

\section{Linearly Preconditioned Methods for Convex Quadratic Objective Functions}
\label{sect:lin:precond}

In this section we discuss the use of linearly preconditioned iterations for the minimization of the convex quadratic objective function \eqref{eq:quadratic:objective}. The optimality equations of this problem are given by
$$
\mtx{g}(\mtx{x})=\mtx{Ax}-\mtx{b}=\mtx{0}.
$$
We consider optimization methods that solve $\mtx{g}(\mtx{x})=\mtx{0}$ by some form of fixed-point iteration. For example, one of the simplest choices for solving $\mtx{g}(\mtx{x})=\mtx{Ax}-\mtx{b}=\mtx{0}$ is Richardson iteration \cite{Richarson1911, kelley1999iterative}:
$$
\mtx{x}_{k+1}=\mtx{x}_k - \mtx{g}(\mtx{x}_k) = \mtx{x}_k - (\mtx{Ax}_k-\mtx{b}).
$$
This is, in fact, equivalent to steepest descent with unit step length, and converges if $\norm{\mtx{I}-\mtx{A}}<1$.

Two different preconditioning strategies are described in subsections \ref{subsect:precond:sor} and \ref{subsect:precond:cov}: we may either
\begin{inparaenum}[(i)]
\item apply a left\hyp{}preconditioning matrix $\mtx{P}$ to the optimality equations and solve the left\hyp{}preconditioned system $\mtx{P}\mtx{A}\mtx{x} = \mtx{P}\mtx{b}$; or
\item introduce a change of variables $\mtx{x} = \mtx{C}\mtx{z}$ and solve the transformed optimization problem with $\widehat{f}(\mtx{z})=f(\mtx{Cz})$.
\end{inparaenum}
In Sections \ref{subsect:precond:cg}--\ref{subsect:precond:broy} we discuss how these strategies define preconditioned CG, L\hyp{}BFGS, and L\hyp{}Broyden iterations, and we extend them to nonlinear preconditioners in Section \ref{sect:nonlinear:preconditioning}.


\subsection{Linear Left Preconditioning (LP)}
\label{subsect:precond:sor}

Instead of solving the optimality equations $\mtx{g}(\mtx{x})=\mtx{0}$, we can apply Richardson iteration to the left-preconditioned optimality equations
$$
\mtx{P} \mtx{g}(\mtx{x})= \mtx{P} (\mtx{Ax}-\mtx{b})=\mtx{0},
$$
to obtain
\begin{equation}
\mtx{x}_{k+1}=\mtx{x}_k-\mtx{P}(\mtx{Ax}_k-\mtx{b}).
\label{eq:iter:Rich}
\end{equation}
Here $\mtx{P}$ could be chosen to be the matrix from any of the stationary linear iterations commonly used as preconditioners, such as Gauss\hyp{}Seidel (GS); successive over\hyp{}relaxation (SOR); or, since we assume $\mtx{A}$ to be SPD, symmetric GS (SGS) and symmetric SOR (SSOR). Using the matrix splitting $\mtx{A} = \mtx{D} + \mtx{L} + \mtx{U}$, SOR is equivalent to \eqref{eq:iter:Rich} with
$$
\mtx{P} = \omega\left(\mtx{D} + \omega\mtx{L}\right)^{-1},
$$
where $\omega \in (0,2)$ and GS corresponds to the particular choice of $\omega=1$. Similarly, the preconditioner matrix for SSOR is
\begin{equation}
\mtx{P} = \left[\left(\mtx{D}+\omega\mtx{U}^\intercal\right)\frac{1}{\omega(2-\omega)}\mtx{D}^{-1}\left(\mtx{D}+\omega\mtx{U}\right)\right]^{-1},
\label{eq:P:SSOR}
\end{equation}
with $\omega=1$ corresponding to SGS.

In the preconditioned update formula \eqref{eq:iter:Rich}, we take a step in direction $\mtx{P} \mtx{g}(\mtx{x})$ instead of the gradient direction $\mtx{g}(\mtx{x})$, and we can interpret $\mtx{P} \mtx{g}(\mtx{x})$ as the preconditioned gradient direction. In the context of optimization problems, this form of left preconditioning works in general by replacing any occurrence of the gradient $\mtx{g}(\mtx{x})$ in the fixed\hyp{}point iteration (such as Richardson) by the preconditioned gradient $\mtx{P}\mtx{g}(\mtx{x})$ (we are applying the fixed-point iteration to $\mtx{P}\mtx{g}(\mtx{x})=\mtx{0}$ instead of $\mtx{g}(\mtx{x})=\mtx{0}$).

\subsection{Linear Transformation Preconditioning (TP)}
\label{subsect:precond:cov}

By defining a linear change of variables $\mtx{x} = \mtx{Cz}$ for some nonsingular matrix $\mtx{C}$, we may rewrite the optimization problem \eqref{eq:opt} for
\begin{equation}
\widehat{f}(\mtx{z})=f(\mtx{Cz}).
\label{eq:fhat}
\end{equation}
We then apply our original optimization method (e.g., Richardson, CG or L\hyp{}BFGS) to \eqref{eq:fhat}, and transform back to the $\mtx{x}$ variables.
In doing so, we employ
$$
\nabla_{\mtx{z}} \widehat{f}(\mtx{z}) = \nabla_{\mtx{z}} f(\mtx{Cz}) =  \mtx{C}^\intercal \nabla_\mtx{x} f(\mtx{x}),
$$
and observe that in the resulting iteration formula for $\mtx{x}$, matrices $\mtx{C}$ and $\mtx{C}^\intercal$ will appear along with $\nabla_{\mtx{x}} f(\mtx{x})=\mtx{g}(\mtx{x})$.
In particular, any products $\mtx{g}^\intercal(\mtx{x}) \mtx{g}(\mtx{x})$ in the iteration formula will be transformed to $\mtx{g}^\intercal(\mtx{x}) \mtx{C C}^\intercal \mtx{g}(\mtx{x})$.

For the specific example of the convex quadratic minimization problem \eqref{eq:quadratic:objective} and Richardson iteration, the transformed objective and gradient functions are
$$
\widehat{f}(\mtx{z}) = \frac{1}{2}\mtx{z}^\intercal\mtx{C}^\intercal\mtx{A}\mtx{C}\mtx{z} - \left(\mtx{C}^\intercal\mtx{b}\right)^\intercal\mtx{z} \quad \text{and} \quad \widehat{\mtx{g}}(\mtx{z}) = \mtx{C}^\intercal\mtx{A}\mtx{C}\mtx{z} - \mtx{C}^\intercal\mtx{b},
$$
with corresponding iteration
$$
\mtx{z}_{k+1} = \mtx{z}_{k} - \left(\mtx{C}^\intercal\mtx{A}\mtx{C}\mtx{z}_{k} - \mtx{C}^\intercal\mtx{b}\right),
$$
which gives, upon transforming back to $\mtx{x}$,
\begin{equation}
\mtx{x}_{k+1} = \mtx{x}_{k} - \mtx{C}\mtx{C}^\intercal\left(\mtx{A}\mtx{x}_{k} - \mtx{b}\right).
\label{eq:SSOR:x}
\end{equation}
If we call the SPD matrix $\mtx{CC}^\intercal$ the preconditioner matrix $\mtx{P}$ and take it to be the SGS or SSOR matrix, we get, for Richardson, the same result as the LP formula \eqref{eq:iter:Rich}. 

An important observation is that for more elaborate optimization methods such as L\hyp{}BFGS, the LP and TP approaches may give different results. For example, any scalar product $\mtx{g}^\intercal(\mtx{x}) \mtx{g}(\mtx{x})$ in the iteration formula will be transformed to $\mtx{g}^\intercal(\mtx{x}) \mtx{P} \mtx{g}(\mtx{x})$ in the TP approach, whereas it will become $\mtx{g}^\intercal(\mtx{x}) \mtx{P}^\intercal \mtx{P} \mtx{g}(\mtx{x})$ in the LP approach. This difference may appear subtle, and intuitively the TP approach may appear preferable since it is more closely aligned with the original optimization problem, but we will see in the preliminary numerical results for the convex quadratic case, and in the general results after extending the approaches to nonlinear preconditioning, that both approaches may have their merits (corresponding also to the findings for nonlinearly preconditioned NCG in \cite{sterck2014nonlinearly, sterck2016nonlinearly}).

\subsection{Linearly Preconditioned CG}
\label{subsect:precond:cg}

Given that CG with the LP strategy would require the matrix $\mtx{PA}$ to be SPD for defining the weighted norm in which the error is minimized, the LP strategy is inappropriate for CG, and thus we only consider the use of TP for CG. This derivation is well-documented in the literature; see, for instance, \cite{hager2006survey}. Writing the CG iteration in terms of $\mtx{z}$,
\begin{align*}
\mtx{z}_{k+1} & = \mtx{z}_{k} + \widehat{\alpha}_k\widehat{\mtx{d}}_k, \\
\widehat{\mtx{d}}_{k+1} & = -\widehat{\mtx{g}}_{k+1} + \widehat{\beta}_{k}\widehat{\mtx{d}}_k, \quad \widehat{\mtx{d}}_0 = -\widehat{\mtx{g}}_0 ,
\end{align*}
then converting back to $\mtx{x}_k$, we obtain
\begin{equation}
\label{eq:pcg:linear:iter}
\begin{aligned}
\mtx{x}_{k+1} & = \mtx{x}_{k} + \alpha_k\mtx{d}_k, \\
\mtx{d}_{k+1} & = -\mtx{P}\mtx{g}_{k+1} + \widehat{\beta}_{k}\mtx{d}_k, \quad \mtx{d}_0 = -\mtx{P}\mtx{g}_0,
\end{aligned}
\end{equation}
as
$$
\widehat{\alpha}_k = -\frac{\widehat{\mtx{d}}_k^\intercal\widehat{\mtx{g}}_k}{\widehat{\mtx{d}}_k^\intercal\mtx{C}^\intercal\mtx{A}\mtx{C}\widehat{\mtx{d}}_k} = -\frac{\mtx{d}_k^\intercal\mtx{C}^{-\intercal}\mtx{C}^\intercal\mtx{g}_k}{\mtx{d}_k^\intercal\mtx{C}^{-\intercal}\mtx{C}^\intercal\mtx{A}\mtx{C}\mtx{C}^{-1}\mtx{d}_k} = -\frac{\mtx{d}_k^\intercal\mtx{g}_k}{\mtx{d}_k^\intercal\mtx{A}\mtx{d}_k} = \alpha_k,
$$
and
$$
\widehat{\beta}_k = \frac{\widehat{\mtx{g}}_{k+1}^\intercal\widehat{\mtx{g}}_{k+1}}{\widehat{\mtx{g}}_k^\intercal\widehat{\mtx{g}}_k} = \frac{\mtx{g}_{k+1}^\intercal\mtx{P}\mtx{g}_{k+1}}{\mtx{g}_k^\intercal\mtx{P}\mtx{g}_k}.
$$

\subsection{Linearly Preconditioned L\hyp{}BFGS}
\label{subsect:precond:bfgs}

The LP version of preconditioned L\hyp{}BFGS is obtained by the direct replacement of each gradient $\mtx{g}_k$ with the left\hyp{}preconditioned gradient $\mtx{P}\mtx{g}_k$ in the components $\mtx{Y}_k$, $\mtx{D}_k$, $\mtx{R}_k$, and $\gamma_k$ of \eqref{eq:inverse:bfgs:compact}, and in computing the QN direction $\mtx{p}_k$. 
%

To derive TP L\hyp{}BFGS, we write \eqref{eq:inverse:bfgs:compact} for $\widehat{f}(\mtx{z})$ as
$$
\widehat{\mtx{H}}_{k} = \widehat{\mtx{H}}_0^{(k)} +
\begin{bmatrix}
\widehat{\mtx{S}}_k & \widehat{\mtx{H}}_0^{(k)}\widehat{\mtx{Y}}_k
\end{bmatrix}
\begin{bmatrix}
\widehat{\mtx{R}}_k^{-\intercal}(\widehat{\mtx{D}}_k + \widehat{\mtx{Y}}_k^\intercal\widehat{\mtx{H}}_0^{(k)}\widehat{\mtx{Y}}_k)\widehat{\mtx{R}}_k^{-1} & -\widehat{\mtx{R}}_k^{-\intercal} \\
-\widehat{\mtx{R}}_k^{-1} & \mtx{0}
\end{bmatrix}
\begin{bmatrix}
\widehat{\mtx{S}}_k^\intercal \\
\widehat{\mtx{Y}}_k^\intercal\widehat{\mtx{H}}_0^{(k)}
\end{bmatrix},
$$
where $\widehat{\mtx{H}}_0^{(k)}=\widehat{\gamma}_k\mtx{I}$. By examining definitions \eqref{eq:defn:Sk}, \eqref{eq:defn:Yk}, \eqref{eq:defn:Dk}, \eqref{eq:defn:Rk}, and \eqref{eq:defn:gamma}, we obtain the following relationships between original and transformed quantities:
$$
\widehat{\mtx{S}}_k = \mtx{C}^{-1}\mtx{S}_k, \quad \widehat{\mtx{Y}}_k = \mtx{C}^{\intercal}\mtx{Y}_k, \quad \widehat{\mtx{D}}_k = \mtx{D}_k, \quad \widehat{\mtx{R}}_k = \mtx{R}_k,
$$
and
\begin{equation}
\widehat{\gamma}_k = \frac{\mtx{s}_{k-1}^\intercal\mtx{y}_{k-1}}{\mtx{y}_{k-1}^\intercal\mtx{P}\mtx{y}_{k-1}}.
\label{eq:defn:gamma:hat} 
\end{equation}
The QN update equation
$$
\mtx{z}_k = \mtx{z}_{k-1} + \alpha_k\widehat{\mtx{H}}_k\widehat{\mtx{g}}(\mtx{z}_{k-1})
$$
transforms to
$$
\mtx{x}_k = \mtx{x}_{k-1} + \alpha_k\mtx{C}\widehat{\mtx{H}}_k\mtx{C}^\intercal\mtx{g}(\mtx{x}_{k-1}).
$$
Computing $\mtx{H}_k := \mtx{C}\widehat{\mtx{H}}_k\mtx{C}^\intercal$, we have
\begin{equation}
\mtx{H}_k = \widehat{\gamma}_k\mtx{P} + 
\begin{bmatrix}
\mtx{S}_k & \widehat{\gamma}_k\mtx{P}\mtx{Y}_k
\end{bmatrix}
\begin{bmatrix}
\mtx{R}_k^{-\intercal}(\mtx{D}_k + \widehat{\gamma}_k\mtx{Y}_k^\intercal\mtx{P}\mtx{Y}_k)\mtx{R}_k^{-1} & -\mtx{R}_k^{-\intercal} \\
-\mtx{R}_k^{-1} & \mtx{0}
\end{bmatrix}
\begin{bmatrix}
\mtx{S}_k^\intercal \\
\widehat{\gamma}_k\mtx{Y}_k^\intercal\mtx{P}
\end{bmatrix}.
\label{eq:precond:bfgs:update}
\end{equation}
Minimizing the preconditioned objective $\widehat{f}(\mtx{z})$ using L\hyp{}BFGS is equivalent to applying L\hyp{}BFGS to $f(\mtx{x})$ where $\mtx{H}_0^{(k)} = \widehat{\gamma}_k\mtx{P}$, which is essentially the same preconditioning strategy described for BFGS in \cite[\S~10.7]{luenberger2015linear}, except they omit the scaling factor $\widehat{\gamma}_k$.

\subsection{Linearly Preconditioned L\hyp{}Broyden}
\label{subsect:precond:broy}

Similar to L\hyp{}BFGS, the LP L\hyp{}Broyden update is obtained by the replacement of each gradient $\mtx{g}_k$ with the left\hyp{}preconditioned gradient $\mtx{P}\mtx{g}_k$ in the component $\mtx{Y}_k$ of \eqref{eq:inverse:broyden:compact} and in computing the QN direction $\mtx{p}_k$. 

To derive the TP L\hyp{}Broyden update we write \eqref{eq:inverse:broyden:update} in terms of $\mtx{z}$ to obtain:
$$
\widehat{\mtx{A}}_k^{-1} = \left[\widehat{\mtx{A}}_0^{(k)}\right]^{-1} - \left(\left[\widehat{\mtx{A}}_0^{(k)}\right]^{-1}\widehat{\mtx{Y}}_k - \widehat{\mtx{S}}_k\right)\left(\widehat{\mtx{M}}_{k} + \widehat{\mtx{S}}_k^\intercal\left[\widehat{\mtx{A}}_0^{(k)}\right]^{-1}\widehat{\mtx{Y}}_k\right)^{-1}\widehat{\mtx{S}}_k^\intercal\left[\widehat{\mtx{A}}_0^{(k)}\right]^{-1},
$$
with $\left[\widehat{\mtx{A}}_0^{(k)}\right]^{-1} = \widehat{\eta}_k\mtx{I}$,
where we take $\widehat{\eta}_k=\widehat{\gamma}_k$ as in \eqref{eq:defn:gamma:hat}. 
Recalling the definition \eqref{eq:defn:Mk} for $i > j$
$$
(\widehat{\mtx{M}}_k)_{i,j} = -\widehat{\mtx{s}}_{i-1}^\intercal\widehat{\mtx{s}}_{j-1} = -\mtx{s}_{i-1}^\intercal\left(\mtx{C}\mtx{C}^\intercal\right)^{-1}\mtx{s}_{j-1} = \mtx{g}_{i-1}^\intercal\mtx{s}_{j-1},
$$
where the last equality follows from \eqref{eq:SSOR:x}. As before, 
$$
\mtx{x}_k = \mtx{x}_{k-1} + \alpha_k\mtx{C}\widehat{\mtx{A}}_k^{-1}\mtx{C}^\intercal\mtx{g}(\mtx{x}_{k-1}),
$$
thus the inverse matrix update is
\begin{equation}
\mtx{A}_k^{-1} := \mtx{C}\widehat{\mtx{A}}_k^{-1}\mtx{C}^\intercal = 
\widehat{\gamma}_k\mtx{P} - \left(\widehat{\gamma}_k\mtx{P}\mtx{Y}_k - \mtx{S}_k\right)\left(\widehat{\mtx{M}}_{k} + \widehat{\gamma_k}\mtx{S}_k^\intercal\mtx{Y}_k\right)^{-1}\mtx{S}_k^\intercal\widehat{\gamma}_k.
\label{eq:precond:broyden:update}
\end{equation}
Compared to the L\hyp{}BFGS case, this is not a full replacement of $\left[\widehat{\mtx{A}}_0^{(k)}\right]^{-1}$ by $\widehat{\gamma}_k\mtx{P}$: only two of the instances involve $\mtx{P}$, the remaining two only require $\widehat{\gamma}_k$.

\subsection{Numerical Results for linear LP and TP methods for Convex Quadratic Functions}
\label{subsect:numerical:results:quadratic}

\begin{figure}[!ht]
\centering
\includegraphics[width=0.8\linewidth]{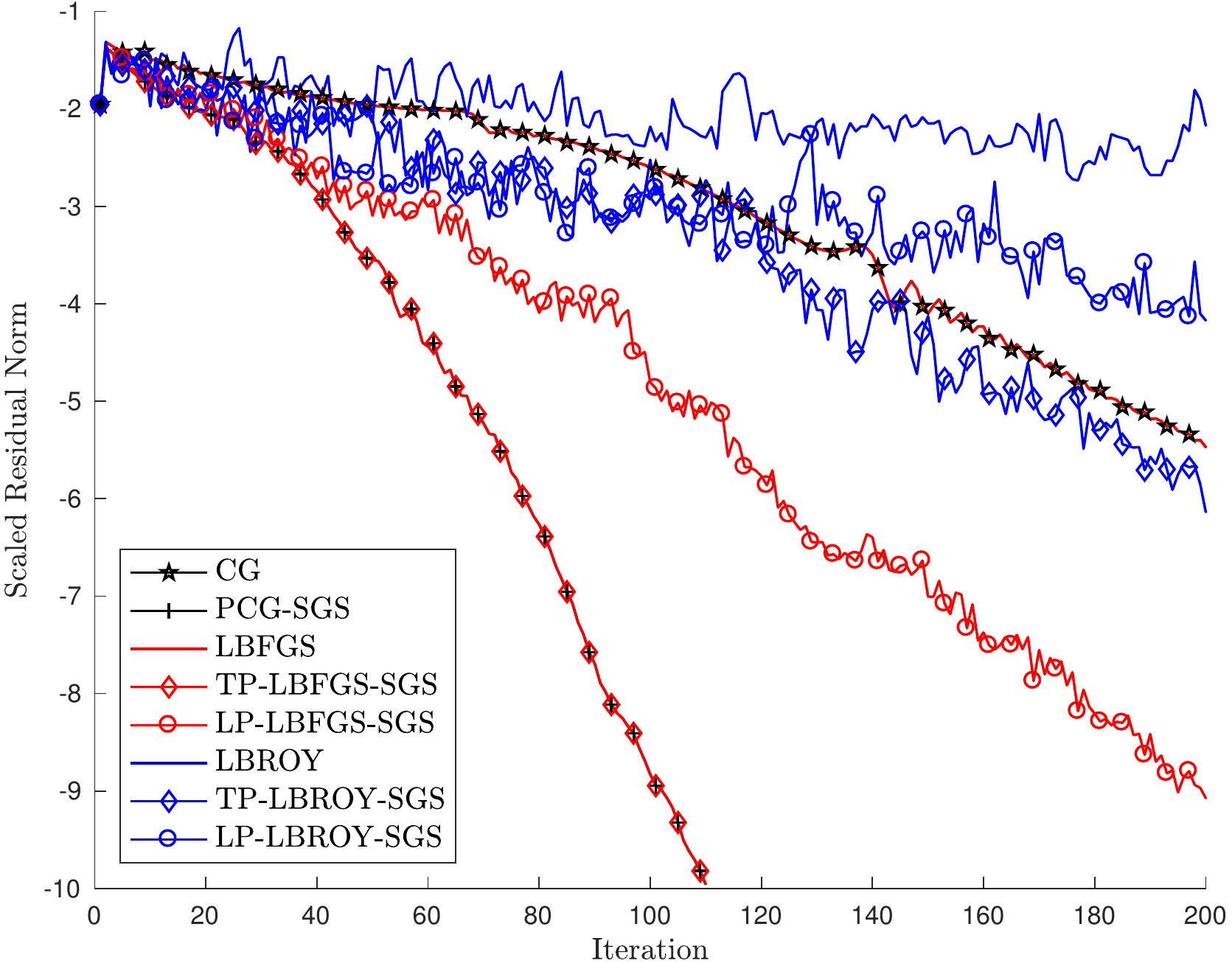}
\caption{Results for CG and QN methods using SGS based preconditioners.}
\label{fig:gs:linear}
\end{figure}

\begin{figure}[!ht]
\centering
\includegraphics[width=0.8\linewidth]{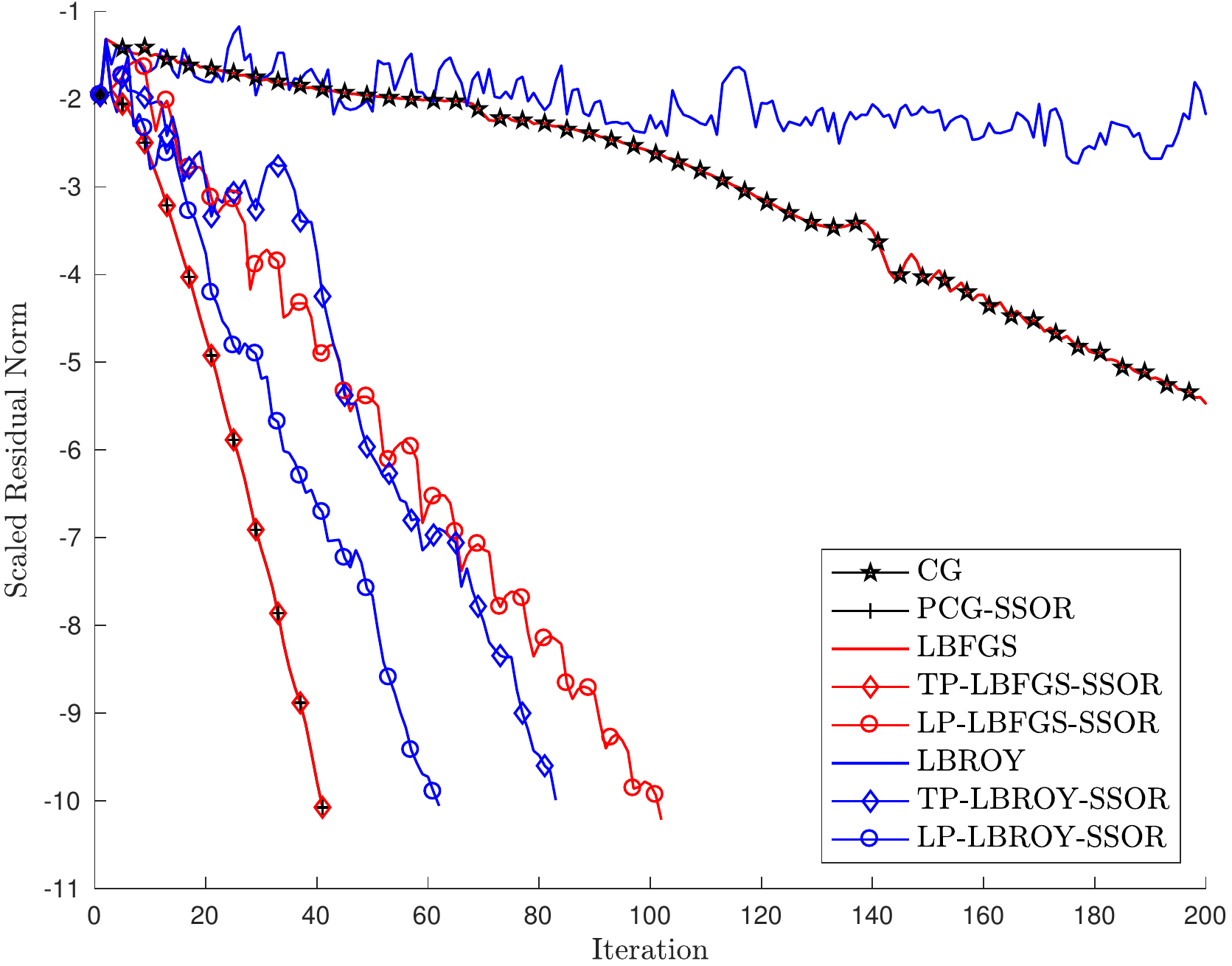}
\caption{Results for CG and QN methods using SSOR based preconditioners.}
\label{fig:sor:linear}
\end{figure}

To illustrate the different preconditioning possibilities, we solve \eqref{eq:quadratic:objective} corresponding to a finite difference discretization of the 2D Poisson equation
$$
u_{xx} + u_{yy} = 2[(1 - 6x^2) y^2 (1 - y^2) + (1 - 6 y^2) x^2 (1 - x^2)], \quad (x,y)\in [0,1]\times[0,1],
$$
with homogeneous Dirichlet boundary conditions and mesh spacing $dx=dy=10^{-2}$, resulting in a problem with $9801$ unknowns. We solve this problem using CG, L\hyp{}BFGS, L\hyp{}Broyden, and their preconditioned variants using SGS or SSOR, the latter with $\omega=1.9$. The QN methods use a window size of $m=5$. For all methods we use the exact step length for quadratic problems $\alpha_k = (\mtx{r}_k^\intercal\mtx{r}_k)/(\mtx{p}_k^\intercal\mtx{A}\mtx{p}_k)$. To precondition L\hyp{}BFGS and L\hyp{}Broyden we consider both LP and TP strategies. These results are presented in Figures \ref{fig:gs:linear} and \ref{fig:sor:linear}, the former containing results for SGS preconditioning and the latter for SSOR preconditioning.

Non\hyp{}preconditioned CG and L\hyp{}BFGS have essentially identical convergence histories, as expected from the discussion of \S~\ref{subsubsect:bfgs:cg}, whereas L\hyp{}Broyden in fact does not converge for this problem, illustrated by the irregular oscillations of the scaled residual norm value. For preconditioned methods, the TP\hyp{}L\hyp{}BFGS overlap the PCG plots for both SGS and SSOR. Larger differences are observed for the LP\hyp{}L\hyp{}BFGS methods, more-so for SGS than SSOR, indicating that preconditioning based on variable transformation is the more effective approach. For L\hyp{}Broyden it is interesting to observe that preconditioning enables the iterations to converge, although the residual curve is still very oscillatory for SGS based preconditioning. When using SGS the TP\hyp{}BROY method is somewhat more effective than the LP version, whereas for SSOR this is reversed, with LP\hyp{}LBROY coming close to the PCG results. Finally, echoing the fact that the SSOR convergence rate is provably better than the convergence rate of SGS \cite{young1954iterative, young1971iterative}, we see that SSOR is clearly a better choice of preconditioner.

\section{Nonlinear Preconditioning Strategies}
\label{sect:nonlinear:preconditioning}


We now consider how to generalize the linear LP and TP strategies to nonlinear preconditioners for general nonlinear optimization problems (\ref{eq:opt}). We first discuss nonlinear preconditioning in general, before describing nonlinearly preconditioned NCG, NGMRES, and QN methods.

\subsection{Nonlinear left preconditioning (LP)}

To generalize the linear LP approach to nonlinear preconditioning, we replace the optimality equation $\mtx{g}(\mtx{x})=\mtx{0}$ with a nonlinearly preconditioned optimality equation
$$
\mathcal{P}(\mtx{g};\mtx{x})=\mtx{0}.
$$
We require solutions of $\mtx{g}(\mtx{x})=\mtx{0}$ to also be solutions of $\mathcal{P}(\mtx{g};\mtx{x})=\mtx{0}$. The notation $\mathcal{P}(\mtx{g};\mtx{x})$ emphasizes that $\mathcal{P}$ is related to solving $\mtx{g}(\mtx{x})=\mtx{0}$. In the convex quadratic case, $\mathcal{P}(\mtx{g};\mtx{x})=\mtx{P} \mtx{g}(\mtx{x}) = \mtx{P}(\mtx{Ax}-\mtx{b})$. In the nonlinear case, $\mathcal{P}(\mtx{g};\mtx{x})$ is generally derived from a nonlinear fixed\hyp{}point equation $\mtx{x} = \mathcal{Q}(\mtx{g};\mtx{x})$, in which case we write
\begin{equation}
\mathcal{P}(\mtx{g};\mtx{x}) := \mtx{x} - \mathcal{Q}(\mtx{g};\mtx{x}).
\label{eq:nonlinear:precond:P}
\end{equation}

Given an iteration $\mtx{x}_{k+1}=\mathcal{Q}(\mtx{g};\mtx{x}_k) = \mtx{x}_k - \mathcal{P}(\mtx{g};\mtx{x}_k)$, we see that
$$
\mathcal{P}(\mtx{g};\mtx{x}_k) = \mtx{x}_k - \mtx{x}_{k+1}
$$
is the (negative) update direction provided by the iteration, which, for a suitable preconditioner should be an improvement on the direction provided by the gradient, $\mtx{g}(\mtx{x})$. In analogy with the linear case where $\mathcal{P}(\mtx{g};\mtx{x}) = \mtx{P}\mtx{g}(\mtx{x})$, we interpret $\mathcal{P}(\mtx{g};\mtx{x})$ as the preconditioned gradient direction. By applying an optimization method with iteration $\mtx{x}_{k+1} = \mathcal{M}(\mtx{g};\mtx{x}_k)$ to $\mathcal{P}(\mtx{g};\mtx{x}) = \mtx{x} - \mathcal{Q}(\mtx{g};\mtx{x}) = \mtx{0}$ instead of to $\mtx{g}(\mtx{x})=\mtx{0}$, we obtain the nonlinearly left\hyp{}preconditioned optimization update
$$
\mtx{x}_{k+1} = \mathcal{M}(\mathcal{P}(\mtx{g};\cdot);\mtx{x}_k).
$$
This means, in practice, that all occurrences of $\mtx{g}(\mtx{x})$ in $\mathcal{M}$ are replaced by $\mathcal{P}(\mtx{g};\mtx{x})$ in the LP approach, as in the case of nonlinear left\hyp{}preconditioning for nonlinear equation systems \cite{brune2015composing}. An important difference in the optimization context, however, is that we continue using the original $f(\mtx{x})$ and $\mtx{g}(\mtx{x})$ in determining the line\hyp{}search step $\alpha$ for methods like CG or L\hyp{}BFGS, so the gradients $\mtx{g}(\mtx{x})$ used in the line\hyp{}search are not replaced by $\mathcal{P}(\mtx{g};\mtx{x})$.

\subsection{Nonlinear transformation preconditioning (TP)}

To extend linear transformation preconditioning to nonlinear preconditioning, consider the iteration formulas derived in Section~\ref{sect:lin:precond} using the linear change of variable $\mtx{x} = \mtx{Cz}$. All occurrences of $\mtx{P}=\mtx{CC}^\intercal$ in the resulting iteration formulas for $\mtx{x}$ appear in front of $\mtx{g}(\mtx{x})$, and we simply replace the linearly preconditioned gradients $\mtx{P} \mtx{g}(\mtx{x})$ by the nonlinearly preconditioned gradients $\mathcal{P}(\mtx{g};\mtx{x})$. This is a natural extension of linear TP to the nonlinear case, with the nonlinear extension reducing to the usual linear preconditioning for nonlinear optimization when $\mathcal{P}(\mtx{g};\mtx{x})$ is chosen as $\mtx{P}\mtx{g}(\mtx{x})=\mtx{CC}^\intercal \mtx{g}(\mtx{x})$ 
\cite{luenberger2015linear,hager2006survey}, and, in the specific case of the NCG method, to the well\hyp{}known formulas for linearly preconditioned CG for SPD linear systems when the objective function is convex quadratic \cite{sterck2014nonlinearly}.

Our approach of nonlinear transformation preconditioning for optimization is based on a change of variables for the solution variable of the nonlinear problem, which is the same general idea behind the nonlinear right-preconditioning method discussed in \cite{brune2015composing} and the related nonlinear elimination preconditioning for the (inexact) Newton method from \cite{LanzkronRoseWilkes1996, HwangSuCai2015}. In the specific context of optimization problems, one could consider the change of variables $\mtx{x}(\mtx{z})$ for the solution variable $\mtx{x}$ in minimizing $f(\mtx{x})$, and then directly minimize $f(\mtx{x}(\mtx{z}))$ for $\mtx{z}$, but this would require computing the Jacobian $\partial \mtx{x} / \partial \mtx{z}$ in $\nabla_{\mtx{z}} f(\mtx{x}(\mtx{z}))=\partial \mtx{x} / \partial \mtx{z} \, \nabla_{\mtx{x}} f(\mtx{x})$, which would in most cases be prohibitively expensive. Instead, in our nonlinear TP approach, we first derive update formulas for the case of a linear change of variables $\mtx{x}(\mtx{z})=\mtx{C}\mtx{z}$, and then, with $\mtx{P}=\mtx{CC}^\intercal$, replace the linearly preconditioned gradients $\mtx{P} \mtx{g}(\mtx{x})$ by the nonlinearly preconditioned gradients $\mathcal{P}(\mtx{g};\mtx{x})$. In this way, we avoid the costly computation of the Jacobian $\partial \mtx{x} / \partial \mtx{z}$. Nonlinear right preconditioning \cite{brune2015composing} and nonlinear elimination preconditioning \cite{LanzkronRoseWilkes1996, HwangSuCai2015} avoid or reduce the extra cost in computing derivatives for the nonlinearly transformed problem in different ways. Also, our nonlinear TP approach results in nonlinearly preconditioned update formulas that, in the case of a linear preconditioning transformation $\mtx{x}=\mtx{C}\mtx{z}$, reduce to well-known linearly preconditioned optimization methods such as PCG for SPD linear systems and linearly preconditioned L\hyp{}BFGS for nonlinear optimization problems, whereas nonlinear right-preconditioning applied to a linear problem results in a new linear system that is different from the system obtained when applying linear right preconditioning \cite{brune2015composing}. The specific details of our nonlinear transformation preconditioning approach for optimization are, thus, different from nonlinear system solvers using nonlinear right-preconditioning as in \cite{brune2015composing} and nonlinear elimination preconditioning as in \cite{LanzkronRoseWilkes1996, HwangSuCai2015}, but these methods are all based on the same idea of modifying the variable of the nonlinear problem.

\subsection{Nonlinearly\hyp{}Preconditioned NCG (NPNCG)}
\label{subsect:left:nonlinear:NCG}

Ideas of using a nonlinear preconditioner with NCG have been around since the 1970s \cite{bartels1974conjugate, concus1978numerical, mittelmann1980on}, but it has not been widely explored. The paper \cite{sterck2014nonlinearly} systematically studied NPNCG iteration in the optimization context, which has since been extended to the matrix manifold setting \cite{sterck2016nonlinearly}. In the nonlinear LP framework (as in \cite{brune2015composing}), given the preconditioning iteration $\mtx{x}_{k+1}=\mathcal{Q}(\mtx{g};\mtx{x}_{k})$, we define
$$
\overline{\mtx{g}}_k := \mathcal{P}(\mtx{g};\mtx{x}_k) = \mtx{x}_k - \mathcal{Q}(\mtx{g};\mtx{x}_k),
$$
and replace every instance of $\mtx{g}_k$ by $\overline{\mtx{g}}_k$ to obtain the LP NPNCG iteration
\begin{equation}
\label{eq:pncg:iter}
\begin{aligned}
\mtx{x}_{k+1} & = \mtx{x}_k + \alpha_k \mtx{p}_k, \\
\mtx{p}_{k+1} & = -\overline{\mtx{g}}_{k+1} + \overline{\beta}_{k+1} \mtx{p}_k, \quad \mtx{p}_0 = -\overline{\mtx{g}}_0.
\end{aligned}
\end{equation}
The corresponding $\beta$ formulas are:
\begin{align}
\label{eq:tildebeta:PR}
\widetilde{\beta}_{k+1}^{\text{PR}} & = \frac{\overline{\mtx{g}}_{k+1}^\intercal\overline{\mtx{y}}_k}{\overline{\mtx{g}}_{k}^\intercal\overline{\mtx{g}}_{k}}, \\
\label{eq:tildebeta:HS}
\widetilde{\beta}_{k+1}^{\text{HS}} & = \frac{\overline{\mtx{g}}_{k+1}^\intercal\overline{\mtx{y}}_k}{\overline{\mtx{y}}_k^\intercal\mtx{p}_{k}}, \\
\label{eq:tildebeta:HZ}
\widetilde{\beta}_{k+1}^{\text{HZ}} & = \left(\overline{\mtx{y}}_k - 2\mtx{p}_k\frac{\norm{\overline{\mtx{y}}_k}^2}{\overline{\mtx{y}}_k^\intercal\mtx{p}_k}\right)^\intercal \frac{\overline{\mtx{g}}_{k+1}}{\overline{\mtx{y}}_k^\intercal\mtx{p}_k},
\end{align}
where $\overline{\mtx{y}}_k = \overline{\mtx{g}}_{k+1} - \overline{\mtx{g}}_{k}$.

To obtain the TP NPNCG iteration, we replace each $\mtx{P}\mtx{g}_k$ in \eqref{eq:pcg:linear:iter} with our preconditioner direction $\overline{\mtx{g}}_{k}$, obtaining \eqref{eq:pncg:iter} again. 
To obtain the TP $\beta$ formulas corresponding to (\ref{eq:beta:PR}--\ref{eq:beta:HZ}), we recall that the products of $\mtx{g}^\intercal\mtx{g}$ transform to $\mtx{g}^\intercal\mtx{CC}^\intercal\mtx{g}$ under the linear transformation $\mtx{x}=\mtx{Cz}$. This becomes $\mtx{g}^\intercal\overline{\mtx{g}}$ when using the nonlinear TP preconditioner, resulting in the following alternatives 
to (\ref{eq:tildebeta:PR}--\ref{eq:tildebeta:HZ})
which incorporate both $\mtx{g}_{k}$ and $\overline{\mtx{g}}_{k}$: 
\begin{align}
\label{eq:hatbeta:PR}
\widehat{\beta}_{k+1}^{\text{PR}} & = \frac{\mtx{g}_{k+1}^\intercal\overline{\mtx{y}}_k}{\mtx{g}_{k}^\intercal\overline{\mtx{g}}_{k}}, \\
\label{eq:hatbeta:HS}
\widehat{\beta}_{k+1}^{\text{HS}} & = \frac{\mtx{g}_{k+1}^\intercal\overline{\mtx{y}}_k}{\mtx{y}_k^\intercal\mtx{p}_{k}}, \\
\label{eq:hatbeta:HZ}
\widehat{\beta}_{k+1}^{\text{HZ}} & = \frac{\mtx{g}_{k+1}^\intercal\overline{\mtx{y}}_{k}}{\mtx{y}_k^\intercal\mtx{p}_k}
- 2\mtx{p}_k^\intercal\mtx{g}_{k+1}
\frac{\mtx{y}_k^\intercal\overline{\mtx{y}}_k}{(\mtx{y}_k^\intercal\mtx{p}_k)^2}.
\end{align}
The LP and TP versions of NPNCG have been considered previously in \cite{sterck2014nonlinearly}. It was observed there numerically that the TP $\widehat{\beta}$ update formulas tend to give better results than the LP $\widetilde{\beta}$ formulas obtained from nonlinear left\hyp{}preconditioning.

\subsection{Nonlinearly Preconditioned L\hyp{}BFGS}
\label{subsect:nonlinear:PQN}

The LP NPQN iteration is obtained by applying the L\hyp{}BFGS update formula to $\mathcal{P}(\mtx{g};\mtx{x}) = 0$. We continue to use \eqref{eq:inverse:bfgs:compact} as our inverse Hessian approximation and define $\overline{\mtx{g}}_k = \mathcal{P}(\mtx{g};\mtx{x}_k)$ and
$\overline{\mtx{y}}_k = \overline{\mtx{g}}_{k+1} - \overline{\mtx{g}}_k$. We replace each instance of $\mtx{g}_k$ with $\overline{\mtx{g}}_k$ and each instance of $\mtx{y}_k$ with $\overline{\mtx{y}}_k$, which preserves the symmetry of $\mtx{H}_{k}$, to obtain 
\begin{equation}
\widetilde{\mtx{H}}_{k} = \widetilde{\gamma}_k + 
\begin{bmatrix}
\mtx{S}_k & \widetilde{\gamma}_k\overline{\mtx{Y}}_k
\end{bmatrix}
\begin{bmatrix}
\overline{\mtx{R}}_k^{-\intercal}(\overline{\mtx{D}}_k + \widetilde{\gamma}_k\overline{\mtx{Y}}_k^\intercal\overline{\mtx{Y}}_k)\overline{\mtx{R}}_k^{-1} & -\overline{\mtx{R}}_k^{-\intercal} \\
-\overline{\mtx{R}}_k^{-1} & \mtx{0}
\end{bmatrix}
\begin{bmatrix}
\mtx{S}_k^\intercal \\
\widetilde{\gamma}_k\overline{\mtx{Y}}_k^\intercal
\end{bmatrix},
\label{eq:nonlinear:left:precond:bfgs:update}
\end{equation}
where
\begin{align}
\overline{\mtx{Y}}_k & = [\overline{\mtx{y}}_{k-m}\,|\, \overline{\mtx{y}}_{k-m+1} \,|\, \cdots \,|\, \overline{\mtx{y}}_{k-1}], \label{eq:defn:Ykbar}\\
\overline{\mtx{D}}_k & = \text{diag}[\mtx{s}_{k-m}^\intercal\overline{\mtx{y}}_{k-m}, \ldots, \mtx{s}_{k-1}^\intercal\overline{\mtx{y}}_{k-1}], \label{eq:defn:Dkbar}\\
(\overline{\mtx{R}}_k)_{i,j} & = \begin{cases}
(\mtx{s}_{k-m-1+i})^\intercal(\overline{\mtx{y}}_{k-m-1+j}) & \text{ if } i \leq j \\
0 & \text{ otherwise }
\end{cases} \label{eq:defn:Rkbar}
\end{align}
and
\begin{equation}
\widetilde{\gamma}_k = \frac{\mtx{s}_{k-1}^\intercal\overline{\mtx{y}}_{k-1}}{\overline{\mtx{y}}_{k-1}^\intercal\overline{\mtx{y}}_{k-1}}.
\label{eq:nonlinear:precond:gamma:tilde}
\end{equation}
If we instead start from the linear TP L\hyp{}BFGS update \eqref{eq:precond:bfgs:update} and replace each $\mtx{P}\mtx{g}_k$ with $\overline{\mtx{g}}_k$ and each $\mtx{P}\mtx{y}_k$ with $\overline{\mtx{y}}_k$, we obtain for the TP NPQN iteration
\begin{equation}
\resizebox{0.9\linewidth}{!}{
$
\widehat{\mtx{H}}_{k}(\mtx{g}_k) = \widehat{\gamma}_k\mathcal{P}(\mtx{g};\mtx{x}_{k}) + 
\begin{bmatrix}
\mtx{S}_k & \widehat{\gamma}_k\overline{\mtx{Y}}_k
\end{bmatrix}
\begin{bmatrix}
\mtx{R}_k^{-\intercal}(\mtx{D}_k + \widehat{\gamma}_k\mtx{Y}_k^\intercal\overline{\mtx{Y}}_k)\mtx{R}_k^{-1} & -\mtx{R}_k^{-\intercal} \\
-\mtx{R}_k^{-1} & \mtx{0}
\end{bmatrix}
\begin{bmatrix}
\mtx{S}_k^\intercal \\
\widehat{\gamma}_k\overline{\mtx{Y}}_k^\intercal
\end{bmatrix}\mtx{g}_k,
$}
\label{eq:nonlinear:precond:bfgs:update}
\end{equation}
where
\begin{equation}
\widehat{\gamma}_k = \frac{\mtx{s}_{k-1}^\intercal\mtx{y}_{k-1}}{\mtx{y}_{k-1}^\intercal\overline{\mtx{y}}_{k-1}}.
\label{eq:nonlinear:precond:gamma}
\end{equation}
The nonlinearly preconditioned quasi\hyp{}Newton search directions are 
$$
\mtx{p}_{k} = -\widetilde{\mtx{H}}_k\overline{\mtx{g}}_k
$$
for LP, and 
$$
\mtx{p}_{k} = -\widehat{\mtx{H}}_k(\mtx{g}_k)
$$
for TP.

Adding nonlinear preconditioning to an L\hyp{}BFGS implementation is relatively straightforward. Once the nonlinear preconditioner $\mathcal{P}(\mtx{g};\mtx{x}_k)$ is defined, the L\hyp{}BFGS QN update step using $\mtx{H}_k$ as in \eqref{eq:inverse:bfgs:compact} is replaced by either of the nonlinearly preconditioned updates \eqref{eq:nonlinear:left:precond:bfgs:update} or \eqref{eq:nonlinear:precond:bfgs:update}. In particular, when using the LP variant \eqref{eq:nonlinear:left:precond:bfgs:update}, we can still make use of the two-loop recursion in Algorithm \ref{alg:2loop}, replacing each $\mtx{g}_k$ and $\mtx{y}_k$ by their nonlinearly preconditioned analogue. This form of left-preconditioned L-BFGS (equivalent to (\ref{eq:nonlinear:left:precond:bfgs:update})) has been considered before in \cite{brune2015composing} in the context of nonlinear systems solvers for PDEs. The nonlinearly transformation-preconditioned form of L-BFGS as in (\ref{eq:nonlinear:precond:bfgs:update}) has, to our knowledge, not been considered before.

Note that, due to the nonlinear preconditioning, the search directions $\mtx{p}_{k}$ are no longer obtained by forming the product of a matrix approximating the inverse Hessian of the objective function, and the gradient $\mtx{g}_k$ of the objective function, as in non-preconditioned L\hyp{}BFGS. The search directions are, instead, formed using the nonlinear expressions $-\widetilde{\mtx{H}}_k\overline{\mtx{g}}_k$ and $-\widehat{\mtx{H}}_k(\mtx{g}_k)$. As such, the nonlinearly preconditioned L\hyp{}BFGS formulas (\ref{eq:nonlinear:left:precond:bfgs:update}) and (\ref{eq:nonlinear:precond:bfgs:update}) can no longer exactly satisfy a secant property and do not maintain the SPD nature of an approximate inverse Hessian, as in non-preconditioned L\hyp{}BFGS. Instead, the resulting nonlinear update formulas derive their motivation from the property that they reduce to well-established linearly preconditioned L\hyp{}BFGS methods in the case of linear preconditioners $\mathcal{P}(\mtx{g};\mtx{x}) = \mtx{P}\mtx{g}(\mtx{x})$. As indicated in our numerical results of Section \ref{sect:numerical:results}, the nonlinearly preconditioned L\hyp{}BFGS methods provide robust and highly efficient optimization methods for difficult tensor problems that substantially outperform the leading existing methods.

\subsection{Nonlinearly Preconditioned L\hyp{}Broyden}

To simplify equations we assume the initial approximation of the inverse Jacobian has the form $[\mtx{A}_0^{(k)}]^{-1} = \eta_k\mtx{I}$, for some scaling factor $\eta_k$. Similar to L\hyp{}BFGS, the nonlinear LP variant is obtained by replacing $\mtx{g}_k$ with $\overline{\mtx{g}}_k = \mathcal{P}(\mtx{g};\mtx{x}_k)$ throughout \eqref{eq:inverse:broyden:compact} and Algorithm \ref{alg:QN}, resulting in
\begin{equation}
\widetilde{\mtx{A}}_k^{-1} = \eta_k\left(\mtx{I} - \left(\eta_k\overline{\mtx{Y}}_k - \mtx{S}_k\right)\left(\mtx{M}_k + \mtx{S}_k^\intercal\eta_k\overline{\mtx{Y}}_k\right)^{-1}\mtx{S}_k^\intercal\right).
\label{eq:nonlinear:left:precond:broyden:update}
\end{equation}
The idea of applying Broyden's method to a fixed point equation has previously been discussed in \cite{fang2009two}, though in the context of nonlinear systems of equations rather than optimization.

If we instead take the linear TP L\hyp{}Broyden update \eqref{eq:precond:broyden:update} into consideration and replace $\mtx{P}\mtx{g}_k$ with $\overline{\mtx{g}}_k$ and $\mtx{P}\mtx{y}_k$ with $\overline{\mtx{y}}_k$, we obtain the operator
\begin{equation}
\widehat{\mtx{A}}_{k}^{-1}(\mtx{g}_k) = \eta_k\left(\mathcal{P}(\mtx{g}_k;\mtx{x}_{k+1}) - \left(\eta_k\overline{\mtx{Y}}_k - \mtx{S}_k\right)\left(\overline{\mtx{M}}_{k} + \eta_k\mtx{S}_k^\intercal\mtx{Y}_k\right)^{-1}\mtx{S}_k^\intercal\mtx{g}_k\right),
\label{eq:nonlinear:precond:broyden:update}
\end{equation}
where
\begin{equation}
(\overline{\mtx{M}}_k)_{i,j} = 
\begin{cases}
\mtx{g}_{i-1}^\intercal\mtx{s}_{j-1} & \text{ if } i > j \\
0 & \text{ otherwise}
\end{cases}.
\end{equation}

As a final note, both TP formulas \eqref{eq:nonlinear:precond:bfgs:update} and \eqref{eq:nonlinear:precond:broyden:update} combine information from the gradient and preconditioner directions, resulting in some additional storage and computational costs when compared to the LP approaches of \eqref{eq:nonlinear:left:precond:bfgs:update} and \eqref{eq:nonlinear:left:precond:broyden:update}.

\section{Numerical Results}
\label{sect:numerical:results}

To compare the nonlinearly preconditioned L\hyp{}BFGS and L\hyp{}Broyden NPQN methods with existing NPNCG and NPNGMRES algorithms from \cite{sterck2012nonlinear, sterck2014nonlinearly, sterck2016nonlinearly} we consider problems of approximating tensors by computing tensor decompositions. We discuss tensor approximation problems based on two commonly used decompositions: the CP decomposition and the Tucker format HOSVD. Details on tensor decompositions are provided in Appendix \ref{appx:tensors}.

Both of these decompositions can be formulated as optimization problems, and both have an ALS\hyp{}type fixed point iteration for computing optimal points. As discussed in \cite{sterck2012nonlinear, sterck2014nonlinearly} for CP tensors and \cite{sterck2016nonlinearly} for Tucker tensors, these fixed point iterations can serve as effective preconditioners $\mathcal{P}(\mtx{g};\mtx{x})$ for NPCG and NPNGMRES, and we shall show this is also the case for NPQN methods.

All of the following experiments were implemented on a MacBook Pro (2.5 GHz Intel Core i7\hyp{}4770HQ, 16 GB 1600 MHz DDR3 RAM) using MATLAB R2016b with the Tensor Toolbox (V2.6) \cite{TTB_Software,TTB_Dense} to handle tensor computations. 

\subsection{CP Decomposition}
\label{subsect:numerical:results:CP}

The rank\hyp{}$R$ CP decomposition of $\mathscr{X} \in \mathbb{R}^{I_1 \times \cdots \times I_N}$ is
\begin{equation}
\llbracket \mtx{A}^{(1)},\mtx{A}^{(2)},\ldots,\mtx{A}^{(N)} \rrbracket \equiv \sum_{r=1}^R \mtx{a}_r^{(1)} \circ \mtx{a}_r^{(2)} \circ \cdots \circ \mtx{a}_r^{(N)},
\label{eq:CP:decomposition}
\end{equation}
where $\mtx{A}^{(n)}\in\mathbb{R}^{I_n\times R}$ for $n=1,\ldots,N$. To compute a CP decomposition we solve
\begin{equation}
\min_{\{\mtx{A}^{(n)}\}} \quad \tfrac{1}{2}\norm{\mathscr{X} - \llbracket \mtx{A}^{(1)},\mtx{A}^{(2)},\ldots,\mtx{A}^{(N)} \rrbracket}_F^2.
\end{equation}
The standard approach is an alternating least squares (ALS) type iteration \cite{carroll1970analysis, harshman1970foundations}. One iteration of CP\hyp{}ALS consecutively updates each of the factor matrices $\mtx{A}^{(i)}$ with the solution of a least-squares problem while keeping the other factors fixed, as summarized in Algorithm \ref{alg:CPALS}.

\begin{algorithm}[ht]
\caption{CP\hyp{}ALS}\label{alg:CPALS}
\begin{algorithmic}[1]
\Procedure{CP\hyp{}ALS}{$\mathscr{X},\mtx{A}^{(1)},\ldots,\mtx{A}^{(N)}$}
\For{$n=1,\ldots,N$} \label{line:CPALS:loop:start}
\State $\pmb{\Gamma}^{(n)} = (\mtx{A}^{(1)\intercal}\mtx{A}^{(1)})\ast \cdots \ast (\mtx{A}^{(n-1)\intercal}\mtx{A}^{(n-1)}) \ast (\mtx{A}^{(n+1)\intercal}\mtx{A}^{(n+1)}) \ast \cdots \ast (\mtx{A}^{(N)\intercal}\mtx{A}^{(N)})$
\State $\mtx{A}^{(n)} = \mtx{X}_{(n)}\left(\mtx{A}^{(N)} \odot \cdots \odot \mtx{A}^{(n+1)} \odot \mtx{A}^{(n-1)} \odot \cdots \odot \mtx{A}^{(1)} \right)^\intercal\left(\pmb{\Gamma}^{(n)}\right)^\dagger$
\EndFor \label{line:CPALS:loop:end}
\State return $\mtx{A}^{(1)},\ldots,\mtx{A}^{(N)}$
\EndProcedure
\end{algorithmic}
\end{algorithm}

\subsubsection{CP Decomposition Results}

We consider standard L\hyp{}BFGS and L\hyp{}Broyden, their left preconditioning (LP) variants \eqref{eq:nonlinear:left:precond:bfgs:update} and \eqref{eq:nonlinear:left:precond:broyden:update}, and their variable transformation preconditioning (TP) variants \eqref{eq:nonlinear:precond:bfgs:update} and \eqref{eq:nonlinear:precond:broyden:update}. For the preconditioner $\mathcal{Q}$ we use either one forward sweep (F) or occasionally one forward\hyp{}backward sweep (FB) of CP\hyp{}ALS (Algorithm \ref{alg:CPALS}, lines~\ref{line:CPALS:loop:start}--\ref{line:CPALS:loop:end} with either $n=1,\ldots,N$ or $n=1,2\ldots,N-1,N,N-1,\ldots,2,1$). 
The forward-backward sweeps of ALS are inspired by SGS and SSOR for the convex quadratic functional, which are forward-backward versions of GS and SOR, respectively.
For the CP decomposition problem we set $\eta_k = 1$ for preconditioned L\hyp{}Broyden methods and $\eta_k = \gamma_k$ as prescribed in \eqref{eq:defn:gamma} for non\hyp{}preconditioned L\hyp{}Broyden, as these were observed to give the best results.

Before comparing to the existing methods for CP decompositions, we first test to determine the best window size $m\in\{1,\ldots,10\}$ and line\hyp{}search method for subsequent experiments. Two line\hyp{}searches are considered. The first is the Mor\'e-Thuente (MT) algorithm from the Poblano Toolbox (v1.1) \cite{dunlavy2010poblano,more1994line}. As in \cite{sterck2012nonlinear, sterck2014nonlinearly}, we use the default parameters: $10^{-4}$ for the sufficient decrease condition tolerance, $10^{-2}$ for the curvature condition tolerance, an initial step length of 1, and a maximum of 20 iterations. 

The second, which we refer to as modified backtracking (modBT), is an attempt at imposing a ``relaxed'' line\hyp{}search condition on the QN step, based on the observation that in certain cases the NPQN methods (and also some of the other methods we compare with) converged faster with a fixed unit step length instead of a line\hyp{}search. In such cases the sequence of objective values was not monotonic, hence modBT does not require the objective value to decrease at every iteration, but only not to increase too much, with the increase tolerated decreasing as iteration count increases. Step lengths of 1, $1/2$, and $1/4$ are considered, accepting the step as soon as growth is small enough; and if all three are rejected, it takes a step of length $1/8$ in the preconditioner direction. This is a feasible approach if methods work with unit step length close to convergence, such as QN, but not for those requiring an accurate line\hyp{}search, such as NCG. This algorithm is summarized in Algorithm \ref{alg:modBT:CP}. The while statement condition assumes our objective values will be non\hyp{}negative, which is the case for the CP decomposition problem. 

\begin{algorithm}[ht]
\caption{Modified Backtracking Line\hyp{}Search}\label{alg:modBT:CP}
\begin{algorithmic}[1]
\Procedure{modBT}{$\mtx{x}_k$,$f_k$,$\mtx{g}_k$,$\mtx{p}_k$,$\text{iter}$}
\State $\alpha_k = 1$, $\text{flag} = 0$
\State $\mtx{x}_{k+1} = \mtx{x}_k + \alpha_k \mtx{p}_k$
\State $f_{k+1} = f(\mtx{x}_{k+1})$
\While{$f_{k+1} > (1+e^{-2\cdot\text{iter}})f_{k}$ \&\& $\text{flag} < 4$}
\State $\text{flag} = \text{flag} + 1$
\If{$\text{flag} == 3$}
\State $\mtx{S}_k = [\,], \mtx{Y}_k = [\,]$
\State $\mtx{p}_k = -\overline{\mtx{g}}_k$
\Else
\State $\alpha_k = 0.5\alpha_k$
\EndIf
\State $\mtx{x}_{k+1} \gets \mtx{x}_k + \alpha_k \mtx{p}_k$
\State $f_{k+1} \gets f(\mtx{x}_{k+1})$
\EndWhile
\State return $\mtx{x}_{k+1}, f_{k+1}, \alpha_k$
\EndProcedure
\end{algorithmic}
\end{algorithm}

Both line\hyp{}searches have a reset condition to recover from bad steps. If the search fails to produce an acceptable step length the QN approximation is reset by clearing $\mtx{S}_k$ and $\mtx{Y}_k$, following which we either take a step in the preconditioner direction or, in the case of the not preconditioned methods, take a step in the steepest descent direction. For MT we repeat the search in this new direction, whereas for modBT we simply take a step of length $\nicefrac{1}{8}$.

To compare these algorithms we compute CP decompositions of an order-$3$ tensor, as in \cite{sterck2012nonlinear, tomasi2006comparison, acar2011scalable}, which is a standard test problem. We form a pseudo\hyp{}random test tensor of size $(I\times I\times I)$ for $I=100$, with known rank $R=5$ and specify the collinearity $C$ of the factors in each mode to be $0.9$, meaning that
\begin{equation}
\frac{\mtx{a}_r^{(n)\intercal}\mtx{a}_s^{(n)}}{\norm{\mtx{a}_r^{(n)}}\norm{\mtx{a}_s^{(n)}}} = C
\end{equation}
for $r\neq s$, $r,s=1,\ldots,R$, and $n=1,2,3$. Highly collinear columns in the factor matrices indicates an ill-conditioned problem, with slow convergence for ALS and other methods.

The methodology for creating such a tensor is described in \cite{tomasi2006comparison}. To this tensor of known rank we add homoskedastic noise (noise with constant variance) and heteroskedastic noise (noise with nonconstant variance). As in \cite{sterck2014nonlinearly, sterck2016nonlinearly}, given $\mathscr{N}_1$ and $\mathscr{N}_2$ with entries from the standard normal distribution, homoskedastic and heteroskedastic noise are added by
\begin{equation}
\mathscr{X}^\prime = \mathscr{X} + \sqrt{\frac{\ell_1}{100-\ell_1}}\frac{\norm{\mathscr{X}}_F}{\norm{\mathscr{N}_1}_F}\mathscr{N}_1 \,\text{ and }\,
\mathscr{X}^{\prime\prime} = \mathscr{X}^\prime + \sqrt{\frac{\ell_2}{100-\ell_2}}\frac{\norm{\mathscr{X^\prime}}_F}{\norm{\mathscr{N}_2 \ast \mathscr{X}^\prime}_F}\mathscr{N}_2 \ast \mathscr{X}^\prime
\label{eq:noise:equations}
\end{equation}
respectively. Parameters $\ell_1$ and $\ell_2$ control noise levels: $\ell_i=0$ corresponding to no noise and $\ell_i=50$ corresponding to noise of the same magnitude as $\mathscr{X}$. For this test we take $\ell_1 = 10$ and $\ell_2 = 1$.

\begin{table}[ht]
\centering
\resizebox{\linewidth}{!}{
\begin{tabular}{lrr|rr|rr|rr|rr|rr|}
\hline
\multicolumn{1}{|l|}{\multirow{1}{*}{Line\hyp{}Search}} & \multicolumn{6}{c|}{MT} & \multicolumn{6}{c|}{modBT} \\ 
\hline
\multicolumn{1}{|l|}{Window Size ($m$)} & \multicolumn{2}{c}{1} & \multicolumn{2}{c}{2} & \multicolumn{2}{c|}{10} & \multicolumn{2}{c}{1} & \multicolumn{2}{c}{2} & \multicolumn{2}{c|}{10} \\
\hline
\multicolumn{1}{l|}{} & Time & Iter & Time & Iter & Time & Iter & Time & Iter & Time & Iter & Time & Iter \\
\hline
\multicolumn{1}{|l|}{L-BFGS} & 5.2 & 941 & 4.9 & 835 & 4.0 & 664 & *3.1 & 1000 & *3.1 & 1000 & *3.2 & 1000 \\
\hline
\multicolumn{1}{|l|}{L-BFGS-LP-F} & 2.2 & 156 & 1.3 & 88 & 1.8 & 112 & \textbf{0.7} & 67 & 0.9 & 85 & \textbf{0.8} & 72 \\
\hline
\multicolumn{1}{|l|}{L-BFGS-TP-F} & \textbf{1.2} & 102 & \textbf{1.1} & 84 & 1.8 & 124 & 0.8 & 79 & 0.8 & 77 & \textbf{0.8} & 72 \\
\hline
\multicolumn{1}{|l|}{L-BROY} & *5.5 & 1000 & *6.8 & 1000 & 6.0 & 933 & *3.3 & 1000 & *4.7 & 1000 & *4.6 & 1000 \\
\hline
\multicolumn{1}{|l|}{L-BROY-LP-F} & 1.8 & 142 & 1.9 & 140 & 2.6 & 183 & \textbf{0.7} & 68 & \textbf{0.7} & 68 & 1.0 & 95 \\
\hline
\multicolumn{1}{|l|}{L-BROY-TP-F} & 1.9 & 164 & 1.5 & 115 & \textbf{1.7} & 131 & \textbf{0.7} & 72 & 0.8 & 83 & 0.9 & 84 \\
\hline
\cline{2-13}
\multicolumn{1}{l|}{} & \multicolumn{12}{c|}{\textbf{Independent of Window Size}} \\ 
\cline{2-13}
\multicolumn{1}{l|}{} & \multicolumn{3}{c|}{Time} & \multicolumn{3}{c|}{Iter} & \multicolumn{3}{c|}{Time} & \multicolumn{3}{c|}{Iter} \\
\hline
\multicolumn{1}{|l|}{CP-ALS} & \multicolumn{3}{c|}{*58.5} & \multicolumn{3}{c|}{10000} & \multicolumn{3}{c|}{*58.5} & \multicolumn{3}{c|}{10000} \\
\hline
\multicolumn{1}{|l|}{NPNGMRES} & \multicolumn{3}{c|}{\textbf{1.6}} & \multicolumn{3}{c|}{100} & \multicolumn{3}{c|}{8.9} & \multicolumn{3}{c|}{452} \\
\hline
\multicolumn{1}{|l|}{NPNCG $\widehat{\beta}_{\text{HS}}$} & \multicolumn{3}{c|}{3.2} & \multicolumn{3}{c|}{144} & \multicolumn{3}{c|}{\textbf{5.0}} & \multicolumn{3}{c|}{612} \\
\cline{1-13}
\multicolumn{1}{|l|}{NPNCG $\widetilde{\beta}_{\text{HS}}$} & \multicolumn{3}{c|}{2.1} & \multicolumn{3}{c|}{107} & \multicolumn{3}{c|}{5.9} & \multicolumn{3}{c|}{708} \\
\hline
\end{tabular}
}
\caption{Results comparing MT and modBT line\hyp{}searches and different window sizes for a small, low-noise CP decomposition problem. Asterisks denote runs which failed to converge. Entries in bold indicate the best time for a given value of $m$.}
\label{table:CP:test:m:cvsrch}
\end{table}

For each combination of window size $m$, solver, and line\hyp{}search we ran ten trials, each corresponding to a different random initial guess for the same random tensor, recording the mean time\hyp{}to\hyp{}solution and number of iterations required. The same set of ten initial guesses were used for all test combinations. The iterations ran until a maximum of $1,\!000$ iterations was reached, $10,\!000$ function evaluations had been computed, or $\norm{\mtx{g}_k}/\numel(\mtx{x})$ decreased below a tolerance of $10^{-7}$, where $\numel(\mtx{x}) = 3IR$ is the number of unknowns in our decomposition. Results for this test are recorded in Table \ref{table:CP:test:m:cvsrch} for MT and modBT. Entries in bold indicate the lowest time for a given value of $m$, and entries with asterisks denote a method that failed to converge to the stated tolerance.

First, it is clear that standard L\hyp{}BFGS and L\hyp{}Broyden fail to converge within the limits imposed in the vast majority of cases. This agrees with previous experiments which observed that L\hyp{}BFGS and NCG without preconditioning do not improve upon ALS in terms of time\hyp{}to\hyp{}solution \cite{acar2011scalable}. In comparison, NCG, NGMRES, and L\hyp{}BFGS methods nonlinearly preconditioned by the ALS iteration perform much better when accurate solutions are desired.

A general trend observed is that NPQN methods using MT tend to require more time to converge compared to those using modBT, whereas the existing NPNCG and NPNGMRES iterations perform better with the MT line\hyp{}search. Increasing window size $m$ may also result in a small increase in computation time. Based on these observations, we will restrict further consideration to NPQN methods using the modBT line\hyp{}search, NPNCG and NPNGMRES using MT, and window sizes of $m=1$ and $m=2$.

Before comparing the selected NPQN methods to NPNCG and NPNGMRES, we provide implementation details for these algorithms, see also \cite{sterck2012nonlinear, sterck2014nonlinearly, sterck2016nonlinearly}. The NPNGMRES least squares system grows until a maximum of $m=20$ past iterates is reached. If NPNGMRES produces an ascent search direction, we restart by discarding all past iterates. NPNCG is restarted by setting $\beta=0$ every 20 iterations. We use the two HS $\beta$ formulas from \eqref{eq:tildebeta:HS} and \eqref{eq:hatbeta:HS}. Successful termination occurs when $\norm{\mtx{g}_k}_F/\numel(\mtx{x}) < 10^{-7}$. Note that NPNCG convergence stalls when $\norm{\mtx{g}}_F \approx 10^{-7}$, a well-known phenomenon for NCG that can be explained by a loss of accuracy in the linesearch step, where the Wolfe sufficient decrease condition is checked \cite{sterck2013steepest, sterck2014nonlinearly, hager2005new}.

\begin{table}[!ht]
\centering
\resizebox{0.4\linewidth}{!}{
\begin{tabular}{c|l|r|r|}
\cline{2-4}
\multicolumn{1}{c}{} & \multicolumn{1}{|l|}{Algorithm \Tstrut} & Time & Iter \\ 
\cline{2-4}
\multicolumn{1}{c}{} & \multicolumn{1}{|l|}{CP-ALS \Tstrut} & *123.9 & 10000 \\
\cline{2-4}
\multicolumn{1}{c}{} & \multicolumn{1}{|l|}{NPNCG $\widehat{\beta}_{\text{HS}}$ \Tstrut} & 5.5 & 91 \\
\cline{2-4}
\multicolumn{1}{c}{} & \multicolumn{1}{|l|}{NPNCG $\widetilde{\beta}_{\text{HS}}$ \Tstrut} & *26.5 & 212 \\
\cline{2-4}
\multicolumn{1}{c}{} & \multicolumn{1}{|l|}{NPNGMRES \Tstrut} & 4.3 & 80 \\
\hline
\multicolumn{1}{|l}{\multirow{7}{*}{$m=1$}} & \multicolumn{1}{|l|}{L-BFGS \Tstrut} & *11.8 & 1000 \\
\cline{2-4}
\multicolumn{1}{|c}{} & \multicolumn{1}{|l|}{L-BFGS-LP-F \Tstrut} & 2.0 & 68 \\
\cline{2-4}
\multicolumn{1}{|c}{} & \multicolumn{1}{|l|}{L-BFGS-TP-F \Tstrut} & 1.9 & 68 \\
\cline{2-4}
\multicolumn{1}{|l}{} & \multicolumn{1}{|l|}{L-BROY \Tstrut} & *1.5 & 900 \\
\cline{2-4}
\multicolumn{1}{|c}{} & \multicolumn{1}{|l|}{L-BROY-LP-F \Tstrut} & 1.7 & 58 \\
\cline{2-4}
\multicolumn{1}{|c}{} & \multicolumn{1}{|l|}{L-BROY-TP-F \Tstrut} & \textbf{1.6} & 59 \\
\hline
\multicolumn{1}{|l}{\multirow{7}{*}{$m=2$}} & \multicolumn{1}{|l|}{L-BFGS \Tstrut} & *11.7 & 1000 \\
\cline{2-4}
\multicolumn{1}{|c}{} & \multicolumn{1}{|l|}{L-BFGS-LP-F \Tstrut} & 2.3 & 72 \\
\cline{2-4}
\multicolumn{1}{|c}{} & \multicolumn{1}{|l|}{L-BFGS-TP-F \Tstrut} & 1.9 & 60 \\
\cline{2-4}
\multicolumn{1}{|l}{} & \multicolumn{1}{|l|}{L-BROY \Tstrut} & *16.7 & 917 \\
\cline{2-4}
\multicolumn{1}{|c}{} & \multicolumn{1}{|l|}{L-BROY-LP-F \Tstrut} & 1.8 & 60 \\
\cline{2-4}
\multicolumn{1}{|c}{} & \multicolumn{1}{|l|}{L-BROY-TP-F \Tstrut} & 2.4 & 74 \\
\hline
\end{tabular}
}
\caption{Results for computing CP decomposition of a synthetic $(200\times200\times200)$ rank\hyp{}$5$ tensor with high collinearity $(C=0.9)$ and with significant heteroskedastic ($\ell_1=20$) and homoskedastic ($\ell_2=10$) noise.  Results correspond to the average of 10 trials. Asterisks denote algorithms which failed to converge at least once. The bold entry indicates the best average time out of all methods tested.}
\label{table:CP:test:full}
\end{table}

As a basis of comparison we again use the previous method to form a test tensor with specified collinearity and noise levels. We decompose an order\hyp{}$3$ tensor of size $(I\times I\times I)$ for $I=200$ with known rank $R=5$, collinearity $C=0.9$, and set noise parameters $\ell_1 = 20$ and $\ell_2 = 10$. For this problem we run ten trials, with each trial corresponding to a different random initial guess. The results for this test are recorded in Table \ref{table:CP:test:full}. The best time for each trial is indicated in bold. 

For this problem we observe that the CP\hyp{}ALS algorithm fails to converge to the desired tolerance for each trial, and that NPNCG fails twice for $\widetilde{\beta}$. Of the pre\hyp{}existing methods NPNGMRES converged consistently and in most cases exhibited the quickest time\hyp{}to\hyp{}solution. For the newly proposed methods, we observed that all NPQN methods tested outperformed the pre\hyp{}existing methods in terms of solution time in nearly all cases, with improvements of more than 50\% for the best new method. The results indicate increasing window size generally increased the time\hyp{}to\hyp{}solution, though there are some exceptions. Overall, the lowest times for the majority of trials are for $m=1$. It is interesting to observe that the nonlinearly preconditioned L\hyp{}Broyden iterations typically gave the best performance in this case, rather than L\hyp{}BFGS.

\subsection{The Tucker HOSVD}
\label{subsect:numerical:results:HOSVD}

A tensor $\mathscr{X}$ is expressed in Tucker format as $(\mtx{A}^{(1)}, \ldots, \mtx{A}^{(N)})\cdot \mathscr{S}$, where $\mathscr{S}$ is a smaller tensor. If we further require that each $\mtx{A}^{(n)}\in\mathbb{R}^{I_n \times R_n}$ is orthogonal we obtain a Tucker HOSVD. The best approximate HOSVD of a given $\mathscr{X}$ can be determined by solving
\begin{equation}
\label{eq:hosvd:max}
\begin{aligned}
\min_{\{\mtx{A}^{(n)}\}} & \quad -\tfrac{1}{2}\norm{\left(\mtx{A}^{(1)\intercal},\ldots,\mtx{A}^{(N)\intercal}\right)\cdot\mathscr{X}}_F^2 \\
\text{subject to} & \quad \mtx{A}^{(n)}\in\mathbb{R}^{I_n\times R_n} \text{ and } \mtx{A}^{(n)\intercal}\mtx{A}^{(n)} = \mtx{I}_{R_n},
\end{aligned}
\end{equation}
where $\mathscr{S} = \left(\mtx{A}^{(1)\intercal},\ldots,\mtx{A}^{(N)\intercal}\right)\cdot\mathscr{X}$, see Appendix~\ref{appx:tensors}. The workhorse algorithm for solving \eqref{eq:hosvd:max} is the higher-order orthogonal iteration (HOOI), first proposed in \cite{de2000best}, which alternatingly updates the matrices $\mtx{A}^{(i)}$ and is summarized in Algorithm \ref{alg:HOOI}.

\begin{algorithm}[ht]
\caption{HOOI}\label{alg:HOOI}
\begin{algorithmic}[1]
\Procedure{HOOI}{$\mathscr{X}, \mtx{A}^{(1)}, \ldots, \mtx{A}^{(N)}$}
\For{$n=1,\ldots,N$} \label{line:HOOI:loop:start}
\State $\mathscr{Y} \gets (\mtx{A}^{(1)\intercal},\ldots,\mtx{A}^{(n-1)\intercal}, \mtx{I}, \mtx{A}^{(n+1)\intercal}, \ldots,\mtx{A}^{(N)\intercal})\cdot\mathscr{X}$
\State $\mtx{A}^{(n)} \gets R_n$ leading left singular vectors of the matricization $\mtx{Y}_{(n)}$ of $\mathscr{Y}$ 
\EndFor \label{line:HOOI:loop:end}
\State $\mathscr{S} \gets (\mtx{A}^{(1)\intercal},\ldots, \mtx{A}^{(N)\intercal}) \cdot \mathscr{X}$
\State return $\mathscr{S}, \mtx{A}^{(1)}, \ldots, \mtx{A}^{(N)}$
\EndProcedure
\end{algorithmic}
\end{algorithm}

As mentioned in the introduction, the existing NPNCG and NPNGMRES algorithms of \cite{sterck2016nonlinearly} for the Tucker HOSVD problem involve matrix manifold optimization techniques. A brief discussion of matrix manifolds and a manifold NPQN algorithm are given in Appendix~\ref{appx:manifolds}.

\subsubsection{Tucker HOSVD Results}

We perform three sets of tests for the Tucker HOSVD problem: the first to determine combinations of methods, line\hyp{}searches, and parameters which work best, and the remaining two to compare existing and newly proposed methods for synthetic and real\hyp{}life data tensors of different sizes and noise levels.

We first considered all possible combinations of L\hyp{}BFGS and L\hyp{}Broyden with left preconditioning, transformation preconditioning, or no preconditioning; this time using HOOI (Algorithm \ref{alg:HOOI}, lines~\ref{line:HOOI:loop:start}--\ref{line:HOOI:loop:end} with either $n=1,\ldots,N$ or $n=1,2\ldots,N,N-1,\ldots,1$) as $\mathcal{Q}(\mtx{g};\mtx{x})$. For L\hyp{}Broyden we set $\eta_k$ to use the $\gamma_k$ formula corresponding to the equivalent L\hyp{}BFGS variant because this gave the best results. To narrow down the set of variants, we compared all methods in terms of choice of line\hyp{}search and the L\hyp{}BFGS and L\hyp{}Broyden methods in terms of window parameter $m\in\{1,2\}$. 

We again consider MT and modBT, using the same parameters as for the CP problem. The same reset conditions are used for MT, and for modBT the while condition now uses $f_{k+1} > (1-e^{-2\cdot\text{iter}})f_k$, as objective values are negative for the Tucker HOSVD problem. In both cases the QN approximation is reset by clearing $\mtx{S}_k$ and $\mtx{Y}_k$. The NPNGMRES system grows to a maximum of $m=25$. If NPNGMRES produces an ascent search direction $\mtx{p}_k$, we discard all past iterates and search in the direction $-\mtx{p}_k$. NPNCG methods are restarted every 50 iterations. We again use the two HS $\beta$ update parameters from \eqref{eq:tildebeta:HS} and \eqref{eq:hatbeta:HS}. Successful termination occurs when $\norm{\mtx{g}_k}_F/|f(\mtx{x}_k)| < 10^{-7}$.

The first tests involved decomposing a medium-size order\hyp{}$3$ tensor of size $28\times28\times2500$ into a rank $(14,14,100)$ approximation. This tensor was formed using a subset of the MNIST Database of Handwritten Digits \cite{mnist}, previously used for Tucker decomposition tests in \cite{savas2007handwritten, vannieuwenhoven2012new, sterck2016nonlinearly}, which is a collection of 70,000 images of digits centered in a $28\times28$ pixel image. Our test tensor consisted of $2500$ images of the digit 5 with significant additive noise from a uniform distribution over $[0,1]$: $\mathscr{X}^\prime = \mathscr{X} + 2.5 \frac{\norm{\mathscr{X}}}{\norm{\mathscr{N}}} \mathscr{N}$. 

Ten trials were ran for each combination, with each trial corresponding to a different tensor $\mathscr{N}$. HOSVD truncation (Algorithm~\ref{alg:HOSVD} in Appendix~\ref{appx:tensors}) was used to generate the initial point for each trial, and iterations ran until reaching a maximum of 250 iterations, a total execution time greater than 1500 seconds, or $\norm{\mtx{g}_k}/|f(\mtx{x}_k)| < 10^{-7}$. When recording computation time, we omitted time spent checking the termination condition as less expensive stopping criteria may be used in practice. The mean time\hyp{}to\hyp{}solution and iterations required for each test case are recorded in Table \ref{table:tucker:test:m:digit:cvsrch:No:transport}. Entries in bold denote the lowest time for a given window size, and asterisks indicate runs that did not converge.

\begin{table}[ht]
\centering
\resizebox{0.7\linewidth}{!}{
\begin{tabular}{lrr|rr|rr|rr|}
\hline
\multicolumn{1}{|l|}{\multirow{1}{*}{Line\hyp{}Search}} & \multicolumn{4}{c|}{MT} & \multicolumn{4}{c|}{modBT} \\ 
\hline
\multicolumn{1}{|l|}{Window Size ($m$)} & \multicolumn{2}{c}{1} & \multicolumn{2}{c|}{2} & \multicolumn{2}{c}{1} & \multicolumn{2}{c|}{2} \\
\hline
\cline{2-9}
\multicolumn{1}{l|}{} & \multicolumn{8}{c|}{\textbf{Hessian Update without Vector Transport}} \\ 
\cline{2-9}
\multicolumn{1}{l|}{} & Time & Iter & Time & Iter & Time & Iter & Time & Iter \\
\hline
\multicolumn{1}{|l|}{L-BFGS} & *155.0 & 251 & *193.6 & 251 & *61.7 & 251 & *65.2 & 251 \\
\hline
\multicolumn{1}{|l|}{L-BFGS-LP-F} & 98.7 & 67 & 133.1 & 95 & 31.6 & 50 & 33.2 & 51 \\
\hline
\multicolumn{1}{|l|}{L-BFGS-TP-F} & \textbf{85.1} & 54 & \textbf{93.8} & 50 & \textbf{28.0} & 39 & \textbf{31.0} & 42 \\
\hline
\multicolumn{1}{|l|}{L-BROY} & *250.1 & 251 & *270.3 & 251 & *60.4 & 251 & *62.7 & 251 \\
\hline
\multicolumn{1}{|l|}{L-BROY-LP-F} & 113.0 & 77 & 123.6 & 94 & 45.8 & 71 & 53.8 & 81 \\
\hline
\multicolumn{1}{|l|}{L-BROY-TP-F} & 96.3 & 60 &  100.0 & 51 & 40.5 & 57 & 34.3 & 47 \\
\hline
\cline{2-9}
\multicolumn{1}{l|}{} & \multicolumn{8}{c|}{\textbf{Hessian Update with Vector Transport}} \\ 
\cline{2-9}
\multicolumn{1}{l|}{} & Time & Iter & Time & Iter & Time & Iter & Time & Iter \\
\cline{1-9}
\multicolumn{1}{|l|}{L-BFGS-LP-F} & 117.0 & 67 & 161.5 & 92 & 46.4 & 50 & 52.8 & 52 \\
\cline{1-9}
\multicolumn{1}{|l|}{L-BFGS-TP-F} & 98.3 & 54 & 113.9 & 51 & 44.8 & 39 & 55.7 & 43 \\
\cline{1-9}
\multicolumn{1}{|l|}{L-BROY-LP-F} & 134.1 & 77 & 160.3 & 96 & 63.0 & 72 & 58.6 & 64 \\
\cline{1-9}
\multicolumn{1}{|l|}{L-BROY-TP-F} & 113.9 & 60 & 116.5 & 53 & 62.9 & 55 & 60.5 & 48 \\
\hline
\cline{2-9}
\multicolumn{1}{l|}{} & \multicolumn{8}{c|}{\textbf{Independent of Window Size}} \\ 
\cline{2-9}
\multicolumn{1}{l|}{} & \multicolumn{2}{c|}{Time} & \multicolumn{2}{c|}{Iter} & \multicolumn{2}{c|}{Time} & \multicolumn{2}{c|}{Iter} \\
\hline
\multicolumn{1}{|l|}{HOOI} & \multicolumn{2}{c|}{112.8} & \multicolumn{2}{c|}{210} & \multicolumn{2}{c|}{112.8} & \multicolumn{2}{c|}{210} \\
\hline
\multicolumn{1}{|l|}{NPNGMRES} & \multicolumn{2}{c|}{44.9} & \multicolumn{2}{c|}{35} & \multicolumn{2}{c|}{\textbf{35.8}} & \multicolumn{2}{c|}{33} \\
\hline
\multicolumn{1}{|l|}{NPNCG $\widehat{\beta}_{\text{HS}}$} & \multicolumn{2}{c|}{\textbf{33.0}} & \multicolumn{2}{c|}{37} & \multicolumn{2}{c|}{146.0} & \multicolumn{2}{c|}{210} \\
\cline{1-9}
\multicolumn{1}{|l|}{NPNCG $\widetilde{\beta}_{\text{HS}}$} & \multicolumn{2}{c|}{67.2} & \multicolumn{2}{c|}{72} & \multicolumn{2}{c|}{130.1} & \multicolumn{2}{c|}{210} \\
\hline
\end{tabular}
}
\caption{Tucker decomposition results comparing MT and modBT line\hyp{}searches with varying $m$, for a medium-size tensor with MNIST data. Asterisks denote runs which failed to converge. Entries in bold indicate the best time for a given value of $m$.}
\label{table:tucker:test:m:digit:cvsrch:No:transport}
\end{table}

When working in a manifold framework, we may or may not use a vector transport operation when updating the Hessian approximation between iterations (as explained in Appendix~\ref{appx:manifolds}). These tables contain results for both of these possibilities, from which it is clear that omitting this vector transport step results in faster methods, in particular when the modBT line\hyp{}search is used. Because of this, we exclude the vector transport option from further consideration. Next, note that the non\hyp{}preconditioned methods again failed to converge within the maximum number of steps in every case, which was expected based on the CP results. When comparing the line\hyp{}search methods, MT results in slower convergence for all but the NPNCG algorithms. With respect to window size $m$, we only consider the values of $m=1$ and $m=2$, because other values tested (not shown) gave larger execution times.

\begin{table}[!ht]
\centering
\resizebox{0.4\linewidth}{!}{
\begin{tabular}{c|l|r|r|}
\cline{2-4}
\multicolumn{1}{c}{} & \multicolumn{1}{|l|}{Algorithm \Tstrut} & Time & Iter \\ 
\cline{2-4}
\multicolumn{1}{c}{} & \multicolumn{1}{|l|}{HOOI \Tstrut} & 24.7 & 652 \\
\cline{2-4}
\multicolumn{1}{c}{} & \multicolumn{1}{|l|}{NPNCG $\widehat{\beta}_{\text{HS}}$ \Tstrut} & 9.5 & 100 \\
\cline{2-4}
\multicolumn{1}{c}{} & \multicolumn{1}{|l|}{NPNCG $\widetilde{\beta}_{\text{HS}}$ \Tstrut} & 10.1 & 103 \\
\cline{2-4}
\multicolumn{1}{c}{} & \multicolumn{1}{|l|}{NPNGMRES \Tstrut} & 11.2 & 91 \\
\hline
\multicolumn{1}{|l}{\multirow{7}{*}{$m=1$}} & \multicolumn{1}{|l|}{L-BFGS \Tstrut} & 9.0 & 278 \\
\cline{2-4}
\multicolumn{1}{|c}{} & \multicolumn{1}{|l|}{L-BFGS-LP-F \Tstrut} & 4.4 & 89 \\
\cline{2-4}
\multicolumn{1}{|c}{} & \multicolumn{1}{|l|}{L-BFGS-TP-F \Tstrut} & 5.6 & 80 \\
\cline{2-4}
\multicolumn{1}{|l}{} & \multicolumn{1}{|l|}{L-BROY \Tstrut} & 11.6 & 341 \\
\cline{2-4}
\multicolumn{1}{|c}{} & \multicolumn{1}{|l|}{L-BROY-LP-F \Tstrut} & 7.0 & 131 \\
\cline{2-4}
\multicolumn{1}{|c}{} & \multicolumn{1}{|l|}{L-BROY-TP-F \Tstrut} & 9.6 & 134 \\
\hline
\multicolumn{1}{|l}{\multirow{7}{*}{$m=2$}} & \multicolumn{1}{|l|}{L-BFGS \Tstrut} & 8.3 & 249 \\
\cline{2-4}
\multicolumn{1}{|c}{} & \multicolumn{1}{|l|}{L-BFGS-LP-F \Tstrut} & \textbf{4.3} & 89 \\
\cline{2-4}
\multicolumn{1}{|c}{} & \multicolumn{1}{|l|}{L-BFGS-TP-F \Tstrut} & 6.2 & 87 \\
\cline{2-4}
\multicolumn{1}{|l}{} & \multicolumn{1}{|l|}{L-BROY \Tstrut} & *55.7 & 1581 \\
\cline{2-4}
\multicolumn{1}{|c}{} & \multicolumn{1}{|l|}{L-BROY-LP-F \Tstrut} & 7.1 & 130 \\
\cline{2-4}
\multicolumn{1}{|c}{} & \multicolumn{1}{|l|}{L-BROY-TP-F \Tstrut} & 9.3 & 127 \\
\hline
\end{tabular}
}
\caption{Results for decomposing a noisy synthetic $(120\times120\times120)$ tensor into a rank\hyp{}$(20,20,20)$ Tucker HOSVD approximation. Results correspond to the average of 10 trials. Asterisks denote algorithms which failed to converge at least once. The bold entry indicates the best average time out of all methods tested.}
\label{table:tucker:test:synthetic}
\end{table}

For our second test we computed rank\hyp{}$(20,20,20)$ HOSVDs of rank\hyp{}$(40,40,40)$ synthetic tensors of size $(120,120,120)$, using noise parameters $\ell_1=\ell_2=10$ and noise tensors $\mathscr{N}_1$, $\mathscr{N}_2$ with elements from the standard normal distribution; see \eqref{eq:noise:equations}. Ten trials using different noise tensors $\mathscr{N}_1$ and $\mathscr{N}_2$ were carried out for each method. The results, recorded in Table~\ref{table:tucker:test:synthetic}, indicate that HOOI is the slowest of the pre\hyp{}existing methods, typically followed by NPNGMRES and then NPNCG, with no clear winner between the $\widehat{\beta}$ and $\widetilde{\beta}$ variants. For this problem we again use the forward HOOI sweep as nonlinear preconditioner. The L\hyp{}BFGS results clearly improve upon the existing solvers. The L\hyp{}Broyden iterations are less effective, though many are still competitive with the NPNCG and NPNGMRES results. We see that the fastest results are for nonlinearly left\hyp{}preconditioned L\hyp{}BFGS. More generally, we see improvements of more than 50\% for the best new method.

\begin{figure}[!ht]
\centering
\includegraphics[width=0.75\linewidth]{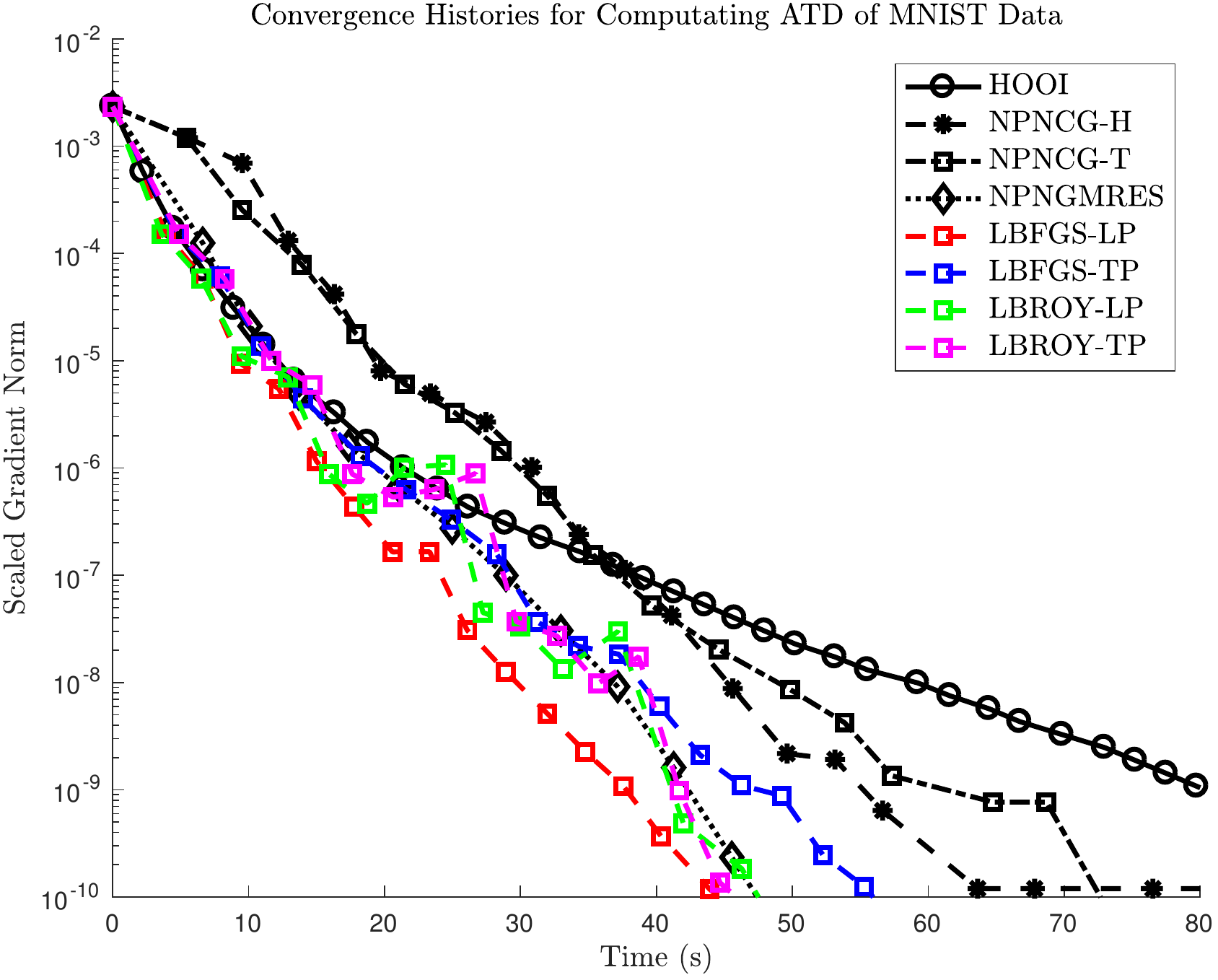}
\includegraphics[width=0.75\linewidth]{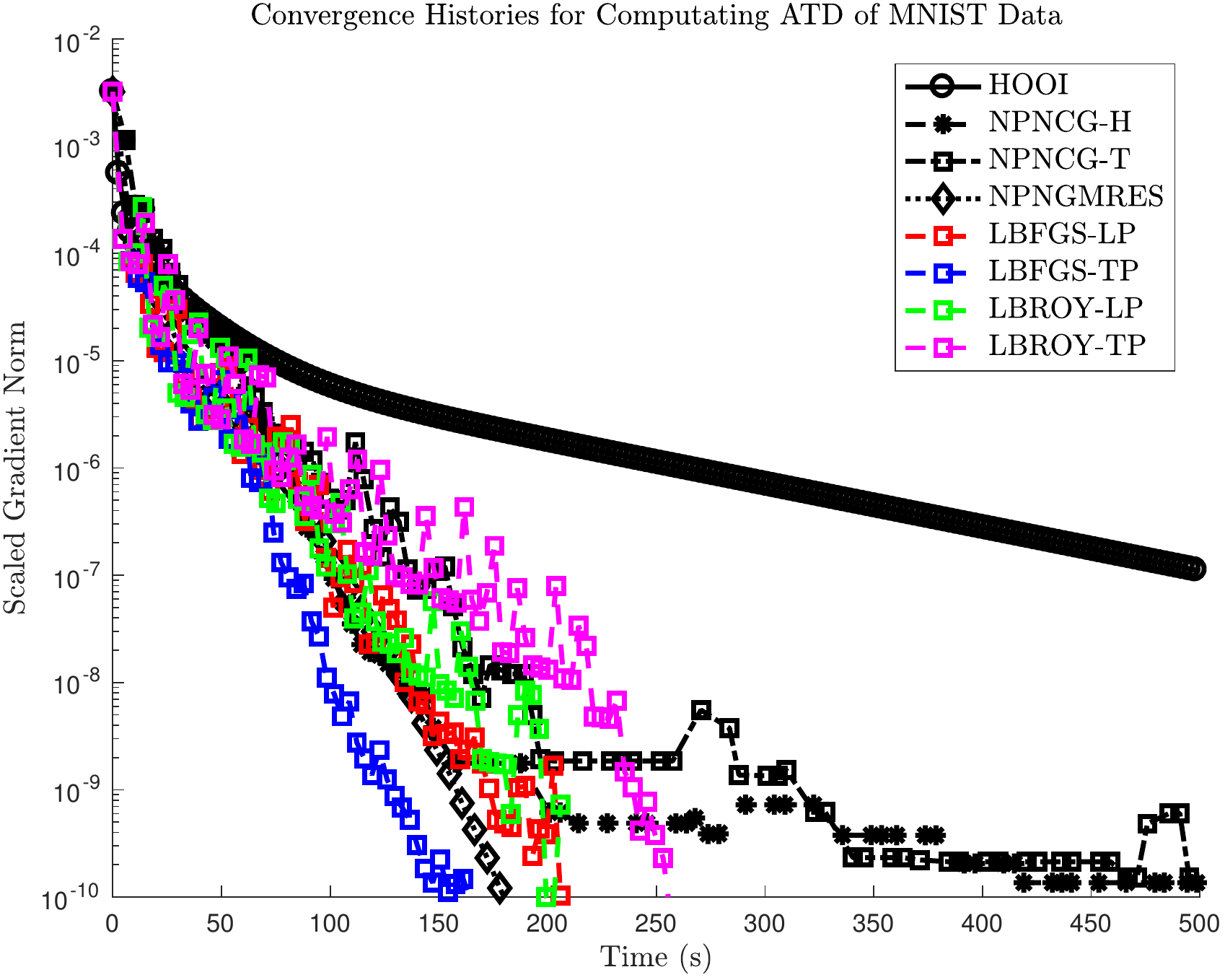}
\caption{Convergence histories showing scaled gradient norms for the rank\hyp{}$(14,14,100)$ Tucker HOSVD decomposition of the $(28\times28\times 5000$) digit tensor without noise (above) and with added noise (below) (corresponding to some of the runs summarized in Table \ref{table:tucker:test:digit}).}
\label{fig:tucker:convergence}
\end{figure}

For the final test we again used the MNIST Database from the initial Tucker test, doubling the number of images to form a tensor $\mathscr{X} \in \mathbb{R}^{28\times28\times 5000}$ consisting of $5000$ images of the digit 5. We add uniformly distributed noise to obtain $\mathscr{X}^\prime = \mathscr{X} + 2.5 \frac{\norm{\mathscr{X}}}{\norm{\mathscr{N}}}\mathscr{N}$, where $\mathscr{N}$ has entries in $[0,1]$. Convergence histories in Figure \ref{fig:tucker:convergence} compare the performance of HOOI, NPNCG, NGMRES, and NPQN methods for a test tensor without (top) and with (bottom) noise. L\hyp{}BFGS and L\hyp{}Broyden without preconditioning are not convergent for these kinds of problems, and hence plots for these solvers are omitted. These plots show that, in the easier noise\hyp{}free case, there is, unsurprisingly, only a small benefit to accelerating HOOI, with NPNGMRES and the NPQN methods all performing slightly better than HOOI, and the NPNCG methods performing slightly worse. Once noise is introduced, however, the convergence of HOOI slows down significantly, and there are clear benefits to using nonlinearly preconditioned methods. In general: nonlinear preconditioning is useful for difficult problems when high accuracy is required, and due to the low amount of overhead involved it does not harm convergence in other circumstances, improving the overall robustness of solvers.

This test is repeated 10 times for different $\mathscr{N}$, with the results recorded in Table \ref{table:tucker:test:digit}. The fastest time is indicated in bold. For existing methods, the general trend is NPNCG using $\widehat{\beta}_{\text{HS}}$ gives the fastest convergence, followed by NPNGMRES, then NPNCG using $\widetilde{\beta}_{\text{HS}}$, and finally HOOI, being the slowest iteration considered by far.
Of the new methods, the L\hyp{}BFGS-TP variant is the fastest, up to 15\% faster than the best previously existing method.
Here we have used forward-backward HOOI sweeps in the nonlinear preconditioners, because for this problem they gave better results than forward sweeps. For L\hyp{}Broyden we note that increasing window size from $m=1$ to $m=2$, results in noticeably faster methods, whereas the difference between mean values for different window sizes is very small for L\hyp{}BFGS.
The 15\% improvement is less than for the previous problem, indicating that this real-data problem may be less ill-conditioned than the artificial-data problem with high collinearity.

\begin{table}[!ht]
\centering
\resizebox{0.4\linewidth}{!}{
\begin{tabular}{c|l|r|r|}
\cline{2-4}
\multicolumn{1}{c}{} & \multicolumn{1}{|l|}{Algorithm \Tstrut} & Time & Iter \\ 
\cline{2-4}
\multicolumn{1}{c}{} & \multicolumn{1}{|l|}{HOOI \Tstrut} & *448.3 & 167 \\
\cline{2-4}
\multicolumn{1}{c}{} & \multicolumn{1}{|l|}{NPNCG $\widehat{\beta}_{\text{HS}}$ \Tstrut} & 115.7 & 29 \\
\cline{2-4}
\multicolumn{1}{c}{} & \multicolumn{1}{|l|}{NPNCG $\widetilde{\beta}_{\text{HS}}$ \Tstrut} & 260.0 & 61 \\
\cline{2-4}
\multicolumn{1}{c}{} & \multicolumn{1}{|l|}{NPNGMRES \Tstrut} & 121.7 & 25 \\
\hline
\multicolumn{1}{|l}{\multirow{7}{*}{$m=1$}} & \multicolumn{1}{|l|}{L-BFGS \Tstrut} & *201.6 & 250 \\
\cline{2-4}
\multicolumn{1}{|c}{} & \multicolumn{1}{|l|}{L-BFGS-LP-FB \Tstrut} & 112.0 & 34 \\
\cline{2-4}
\multicolumn{1}{|c}{} & \multicolumn{1}{|l|}{L-BFGS-TP-FB \Tstrut} & 101.0 & 29 \\
\cline{2-4}
\multicolumn{1}{|l}{} & \multicolumn{1}{|l|}{L-BROY \Tstrut} & *196.0 & 250 \\
\cline{2-4}
\multicolumn{1}{|c}{} & \multicolumn{1}{|l|}{L-BROY-LP-FB \Tstrut} & 139.2 & 43 \\
\cline{2-4}
\multicolumn{1}{|c}{} & \multicolumn{1}{|l|}{L-BROY-TP-FB \Tstrut} & 133.5 & 38 \\
\hline
\multicolumn{1}{|l}{\multirow{7}{*}{$m=2$}} & \multicolumn{1}{|l|}{L-BFGS \Tstrut} & *204.7 & 250 \\
\cline{2-4}
\multicolumn{1}{|c}{} & \multicolumn{1}{|l|}{L-BFGS-LP-FB \Tstrut} & 112.7 & 34 \\
\cline{2-4}
\multicolumn{1}{|c}{} & \multicolumn{1}{|l|}{L-BFGS-TP-FB \Tstrut} & \textbf{99.5} & 28 \\
\cline{2-4}
\multicolumn{1}{|l}{} & \multicolumn{1}{|l|}{L-BROY \Tstrut} & *202.7 & 250 \\
\cline{2-4}
\multicolumn{1}{|c}{} & \multicolumn{1}{|l|}{L-BROY-LP-FB \Tstrut} & 123.0 & 37 \\
\cline{2-4}
\multicolumn{1}{|c}{} & \multicolumn{1}{|l|}{L-BROY-TP-FB \Tstrut} & 109.2 & 31 \\
\hline
\end{tabular}
}
\caption{Results for decomposing a noisy $(28\times28\times 5000)$ tensor into a rank\hyp{}$(100,14,14)$ Tucker HOSVD approximation. Results correspond to the average of 10 trials. Asterisks denote algorithms which failed to converge to the desired tolerance at least once. The bold entry indicates the best average time out of all methods tested.}
\label{table:tucker:test:digit}
\end{table}

\section{Conclusions}
\label{sect:conclusions}

Nonlinear preconditioning strategies are an effective way to improve the convergence of iterative solvers for nonlinear systems and nonlinear optimization problems. In particular, when the problem formulation naturally suggests a fixed\hyp{}point iteration that is more effective than steepest descent, such iterations can be greatly accelerated through use as nonlinear preconditioners. Nonlinearly preconditioned methods can be based on nonlinear left\hyp{}preconditioning or nonlinear preconditioning formulations derived from a change of variables in the optimization problem. 

In this paper we have developed NPQN methods based on the L\hyp{}Broyden and L\hyp{}BFGS update equations. These iterations are applied to the problems of computing two popular tensor decompositions: the CP decomposition and the Tucker HOSVD. These decompositions are commonly used tools in data compression and multilinear statistical analysis, hence improved computational algorithms will continue to be in demand. We also illustrated how to extend these NPQN methods to the manifold setting, which can be useful for solving problems where the unknowns have some constraints, such as underlying symmetry or orthogonal structure.
 
Numerical results provide evidence that, much like the previously developed NPNCG and NPNGMRES methods of \cite{sterck2012nonlinear, sterck2014nonlinearly, sterck2016nonlinearly}, NPQN methods are effectively combined with CP\hyp{}ALS iterations for the CP decomposition and HOOI for the Tucker HOSVD. L\hyp{}BFGS or L\hyp{}Broyden preconditioned by CP\hyp{}ALS or HOOI is always much faster than the individual QN method or ALS\hyp{}type fixed point iteration for difficult problems or when high accuracy is required. Furthermore, our results show that the proposed NPQN methods may significantly outperform NPNCG and NPNGMRES, being up to 50\% faster, establishing them as state of the art methods for difficult, ill\hyp{}conditioned tensor decomposition problems.

There are a number of directions to carry out future work based on these results. At this point we have established the effectiveness of nonlinearly preconditioned versions of the popular NCG, NGMRES, L\hyp{}BFGS and L\hyp{}Broyden methods in Euclidean space and on Grassmann matrix manifolds. These methods can be applied to other systems of equations or optimization problems for matrix or tensor problems that have associated fixed point iterations that are more effective than steepest descent. We can also consider the development of preconditioned versions of other algorithms, such as those based on trust region strategies.

All code and test examples will be made available on the authors' website.

\bibliographystyle{wileyj}
\bibliography{NPQNbib}

\appendix

\section{Tensor Decomposition Problems}
\label{appx:tensors}

\subsection{Tensors}
\label{subsect:tensor:preliminaries}

A \emph{tensor} is a multidimensional array, and the number of \emph{dimensions} (\emph{modes}) of a tensor is called the tensor \emph{order}. Order\hyp{}1 tensors are vectors and order\hyp{}2 tensors are matrices. Tensors of order\hyp{}3 or greater will be indicated by Euler script letters ($\mathscr{X}$) and tensor elements are indicated by subindices or bracketed arguments: $\mathscr{X}_{ijk} = \mathscr{X}(i,j,k)$.

Tensors are useful when large quantities of data need to be organized and analyzed, because each dimension can represent a parameter and each element can represent an observation for a particular parameter combination. As a result tensors have seen widespread use in areas such as in chemometrics \cite{smilde2005multi}, data mining \cite{kolda2008scalable}, food science \cite{bro1998multi}, pattern and image recognition \cite{savas2007handwritten, zhang2013robust}, and signal processing \cite{cichocki2014tensor}.

A \emph{tensor decomposition} expresses a tensor as a sum or product of several low-dimensional components with the goal of simplifying further work involving the tensor data. \emph{Tensor approximation problems} involve seeking the best approximation of a tensor $\mathscr{X}$ by a tensor $\widehat{\mathscr{X}}$, commonly having a specified decomposition, the components of which are determined by minimizing $\norm{\mathscr{X} - \widehat{\mathscr{X}}}$.

\subsubsection{Matrix Singular Value Decomposition (SVD)}
A matrix $\mtx{M}\in \mathbb{R}^{m \times n}$ has SVD $\mtx{M} = \mtx{U\Sigma}\mtx{V}^\intercal$, where $\mtx{U}\in \mathbb{R}^{m \times m}$, $\mtx{V}\in \mathbb{R}^{n \times n}$ are orthogonal and $\mtx{\Sigma}\in \mathbb{R}^{m \times n}$ is diagonal with nonnegative real entries in decreasing order. The nonzero entries of $\mtx{\Sigma}$ are the \emph{singular values} of $\mtx{M}$ and the columns of $\mtx{U}$ $\left(\mtx{V}\right)$ are the \emph{left (right) singular vectors}. The \emph{rank} of $\mtx{M}$ is equal to the number of singular values, and by the Eckhart\hyp{}Young theorem the best rank\hyp{}$r$ approximation of $\mtx{M}$ in the Frobenius norm is obtained by keeping the largest $r$ singular values, setting the rest to zero \cite{eckart1936approximation}.

\subsubsection{Tensor Matricizations and Definitions of Rank}

Mode-$n$ tensor \emph{fibers} are vectors obtained by fixing all indices but the $n^\text{th}$. The \emph{mode\hyp{}$n$ matricization} of $\mathscr{X}$, denoted $\mtx{X}_{(n)}$, has the mode\hyp{}$n$ fibers of $\mathscr{X}$ as its columns. So long as the ordering of fibers is consistent throughout calculations, the specific ordering used is unimportant in many applications. A \emph{rank\hyp{}one tensor} $\mathscr{X}\in\mathbb{R}^{I_1 \times \cdots \times I_N}$ is the outer product
$$
\mathscr{X} = \mtx{a}^{(1)} \circ \mtx{a}^{(2)} \circ \cdots \circ \mtx{a}^{(N)}
$$
where $\mtx{a}^{(n)}\in\mathbb{R}^{I_n}$ for $n=1,\ldots,N$ and
$$
\mathscr{X}_{i_i i_2\ldots i_N} = a_{i_1}^{(1)}a_{i_2}^{(2)}\cdots a_{i_N}^{(N)} \quad \text{for all } 1 \leq i_n \leq I_n.
$$
The \emph{rank} of $\mathscr{X}$, denoted $\tenrank(\mathscr{X})$, is the minimum number of rank\hyp{}one tensors required to express $\mathscr{X}$ as a linear combination \cite{kolda2009tensor}. The \emph{$n$\hyp{}rank} of $\mathscr{X}$ is the dimension of the space spanned by the mode\hyp{}$n$ fibers: $\text{rank}_n(\mathscr{X}) = \text{dim}(\text{Col}(\mtx{X}_{(n)}))$ \cite{de2000multilinear}. The \emph{multilinear rank} of $\mathscr{X}$ is the $N$\hyp{}tuple $(\text{rank}_1(\mathscr{X}),\ldots,\text{rank}_N(\mathscr{X}))$. 

\subsubsection{Tensor and Matrix Products}
\label{subsect:tenprodnorm}

The \emph{mode\hyp{}$n$ contravariant product} of $\mathscr{X}\in\mathbb{R}^{I_1 \times \cdots \times I_N}$ and $\mtx{A}\in\mathbb{R}^{J \times I_n}$ is $\mathscr{Y} = (\mtx{A})_n \cdot \mathscr{X}$ \cite{elden2009newton}:
$$
\mathscr{Y}(i_1,\ldots, i_{n-1}, j, i_{n+1},\ldots, i_n) = \sum_{i_n=1}^{I_n} \mtx{A}(j,i_n) \mathscr{X}(i_1,\ldots, i_N).
$$
Each mode\hyp{}$n$ fiber of $\mathscr{X}$ is multiplied by each row of $\mtx{A}$: $\mtx{Y}_{(n)} = \mtx{AX}_{(n)}.$ It follows that $(\mtx{B})_n \cdot ((\mtx{A})_n \cdot \mathscr{X}) = (\mtx{BA})_n \cdot \mathscr{X}$, and that multiplication in different modes is commutative. 

The \emph{inner product} of $\mathscr{X}$ and $\mathscr{Y}\in\mathbb{R}^{I_1\times\cdots\times I_N}$ is
$$
\langle \mathscr{X}, \mathscr{Y} \rangle = \sum_{i_1=1}^{I_1}\cdots\sum_{i_N=1}^{I_N} \mathscr{X}(i_1, ,\ldots, i_N) \mathscr{Y}(i_1,\ldots, i_N).
$$
The tensor Frobenius norm is $\norm{\mathscr{X}}_F = \sqrt{\langle \mathscr{X}, \mathscr{X} \rangle}$, and $\norm{\mathscr{X}}_F = \norm{\mtx{X}_{(k)}}_F$ for $k=1,\ldots,N$. It is also invariant under orthogonal transformations $\mtx{A}^{(n)}$:
$$
\norm{\mathscr{X}}_F = \norm{(\mtx{A}^{(1)},\ldots,\mtx{A}^{(N)})\cdot\mathscr{X}}_F.
$$

The \emph{Hadamard (element\hyp{}wise) product} of equal sized tensors $\mathscr{X}$ and $\mathscr{Y}$ is denoted $\mathscr{X} \ast \mathscr{Y}$. The \emph{Kronecker} product of $\mtx{A}\in\mathbb{R}^{I\times J}$ and $\mtx{B}\in\mathbb{R}^{K\times L}$ is denoted by $\mtx{A}\otimes\mtx{B}\in\mathbb{R}^{(IK)\times(JL)}$, and the \emph{Khatri\hyp{}Rao product} of $\mtx{C}\in\mathbb{R}^{I\times K}$ and $\mtx{D}\in\mathbb{R}^{J\times K}$ is denoted by $\mtx{C}\odot\mtx{D}\in\mathbb{R}^{(IJ)\times K}$:
$$
\mtx{A}\otimes\mtx{B} = 
\begin{bmatrix}
a_{11}\mtx{B} & a_{12}\mtx{B} & \cdots & a_{1J}\mtx{B} \\
a_{21}\mtx{B} & a_{22}\mtx{B} & \cdots & a_{2J}\mtx{B} \\
\vdots & \vdots & \ddots & \vdots \\
a_{I1}\mtx{B} & a_{I2}\mtx{B} & \cdots & a_{IJ}\mtx{B}
\end{bmatrix},
\quad 
\mtx{C}\odot\mtx{D} = 
\begin{bmatrix}
\mtx{c}_1 \otimes \mtx{d}_1 & \mtx{c}_2 \otimes \mtx{d}_2 & \cdots & \mtx{c}_K \otimes \mtx{d}_K
\end{bmatrix}.
$$
These products are useful when computing tensor decompositions and matricizing tensor\hyp{}matrix products, as shown in the following subsections.

\subsection{CP Decomposition}
\label{subsect:CP:decomposition}

The CP decomposition, also known by the names CANDECOMP (canonical decomposition) and PARAFAC (parallel factors), decomposes a tensor into a sum of rank\hyp{}one tensors \cite{kolda2009tensor}. The rank\hyp{}$R$ CP decomposition of $\mathscr{X} \in \mathbb{R}^{I_1 \times \cdots \times I_N}$ is \cite{kolda2009tensor}
$$
\llbracket \mtx{A}^{(1)},\mtx{A}^{(2)},\ldots,\mtx{A}^{(N)} \rrbracket \equiv \sum_{r=1}^R \mtx{a}_r^{(1)} \circ \mtx{a}_r^{(2)} \circ \cdots \circ \mtx{a}_r^{(N)},
$$
and to compute a CP decomposition we solve
$$
\min_{\{\mtx{A}^{(n)}\}} \quad \tfrac{1}{2}\norm{\mathscr{X} - \llbracket \mtx{A}^{(1)},\mtx{A}^{(2)},\ldots,\mtx{A}^{(N)} \rrbracket}_F^2.
$$
The mode\hyp{}$n$ matricization of a CP tensor is given by
$$
\mtx{A}^{(n)}\left(\mtx{A}^{(N)} \odot \cdots \odot \mtx{A}^{(n+1)} \odot \mtx{A}^{(n-1)} \odot \cdots \odot \mtx{A}^{(1)} \right)^\intercal.
$$
By fixing all matrices except $\mtx{A}^{(n)}$ the problem becomes
\begin{equation}
\min_{\mtx{A}^{(n)}} \quad \tfrac{1}{2}\norm{\mathscr{X}_{(n)} - \mtx{A}^{(n)}(\mtx{A}^{(N)} \odot \cdots \odot \mtx{A}^{(n+1)} \odot \mtx{A}^{(n-1)} \odot \cdots \odot \mtx{A}^{(1)} )^\intercal}_F^2,
\end{equation}
a linear least squares problem with exact solution
$$
\mtx{A}^{(n)} = \mtx{X}_{(n)}\left(\left(\mtx{A}^{(N)} \odot \cdots \odot \mtx{A}^{(n+1)} \odot \mtx{A}^{(n-1)} \odot \cdots \odot \mtx{A}^{(1)} \right)^\intercal\right)^\dagger,
$$
where $\dagger$ denotes the Moore\hyp{}Penrose pseudoinverse \cite{moore1920on, penrose1955generalized}. Since the pseudoinverse of a Khatri\hyp{}Rao product satisfies the identity \cite{kolda2009tensor}
$$
(\mtx{A}\odot\mtx{B})^\dagger = \left(\left(\mtx{A}^\intercal\mtx{A}\right)\ast \left(\mtx{B}^\intercal\mtx{B}\right)\right)(\mtx{A}\odot\mtx{B})^\intercal,
$$
this exact solution is typically implemented as
$$
\mtx{A}^{(n)} = \mtx{X}_{(n)}\left(\mtx{A}^{(N)} \odot \cdots \odot \mtx{A}^{(n+1)} \odot \mtx{A}^{(n-1)} \odot \cdots \odot \mtx{A}^{(1)} \right)^\intercal\left(\pmb{\Gamma}^{(n)}\right)^\dagger,
$$
where
$$
\pmb{\Gamma}^{(n)} = \left(\mtx{A}^{(1)\intercal}\mtx{A}^{(1)}\right)\ast \cdots \ast \left(\mtx{A}^{(n-1)\intercal}\mtx{A}^{(n-1)}\right) \ast \left(\mtx{A}^{(n+1)\intercal}\mtx{A}^{(n+1)}\right) \ast \cdots \ast \left(\mtx{A}^{(N)\intercal}\mtx{A}^{(N)}\right)
$$
for $n=1,\ldots,N$. The CP\hyp{}ALS iteration based on these computations (Algorithm \ref{alg:CPALS}) can be slow to converge in practice, thus alternative optimization algorithms are desirable. Most optimization algorithms require the gradient of
$$
f(\mtx{x}) = \tfrac{1}{2}\norm{\mathscr{X} - \llbracket \mtx{A}^{(1)},\mtx{A}^{(2)},\ldots,\mtx{A}^{(N)} \rrbracket}_F^2,
$$
where $\mtx{x} = (\mtx{A}^{(1)},\mtx{A}^{(2)},\ldots,\mtx{A}^{(N)})$ is the $N$\hyp{}tuple of factor matrices. The partial derivative of $f$ with respect to $\mtx{A}^{(n)}$ is \cite{acar2011scalable}
$$
\frac{\partial f}{\partial \mtx{A}^{(n)}} = -\mtx{X}_{(n)}\left(\mtx{A}^{(N)} \odot \cdots \odot \mtx{A}^{(n+1)} \odot \mtx{A}^{(n-1)} \odot \cdots \odot \mtx{A}^{(1)} \right) + \mtx{A}^{(n)}\pmb{\Gamma}^{(n)}.
$$
Note that setting the gradient of $f$ equal to zero gives
$$
\mtx{A}^{(n)}\pmb{\Gamma}^{(n)} = \mtx{X}_{(n)}\left(\mtx{A}^{(N)} \odot \cdots \odot \mtx{A}^{(n+1)} \odot \mtx{A}^{(n-1)} \odot \cdots \odot \mtx{A}^{(1)} \right),
$$
from which the CP\hyp{}ALS iteration immediately follows.

\subsection{The Tucker HOSVD}
\label{subsect:tucker:HOSVD}

The Tucker format was introduced by Tucker in 1963 for $3$-mode tensors \cite{tucker1963implications}, and has since been extended to $N$-mode tensors; see, for example, \cite{tucker1966some, de2000multilinear}. A tensor $\mathscr{X} \in \mathbb{R}^{I_1 \times \cdots \times I_N}$ is expressed in Tucker format as $(\mtx{A}^{(1)}, \ldots, \mtx{A}^{(N)})\cdot \mathscr{S}$, where $\mathscr{S} \in \mathbb{R}^{R_1 \times \cdots \times R_N}$ and $\mtx{A}^{(n)}\in\mathbb{R}^{I_n \times R_n}$. We must have $R_n \leq I_n$, and in practice often have $R_n \ll I_n$, resulting in a significant reduction in storage. If $R_n \geq \text{rank}_n(\mathscr{X})$ for all $n$, the decomposition is exact. If not, then this is an approximate Tucker decomposition. Tucker decompositions are not unique: replacing $\mathscr{S}$ by $(\mtx{B})_n\cdot\mathscr{S}$ and $\mtx{A}^{(n)}$ by $\mtx{A}^{(n)}\mtx{B}^{-1}$ produces an equivalent tensor.

In \cite{de2000multilinear} the authors introduce a tucker decomposition called the HOSVD and prove all tensors have such a decomposition. The HOSVD of $\mathscr{X}\in \mathbb{R}^{I_1 \times \cdots \times I_N}$ is $\mathscr{X} = (\mtx{A}^{(1)}, \ldots, \mtx{A}^{(N)})\cdot\mathscr{S}$, where $\mathscr{S}\in\mathbb{R}^{I_1 \times \cdots \times I_N}$, each $\mtx{A}^{(n)}\in\mathbb{R}^{I_n \times I_n}$ is orthogonal, and $\mathscr{S}$ satisfies
\begin{enumerate}[(i)]
\item all-orthogonality: for all possible $n$, $\alpha$, and $\beta$, $\alpha\neq\beta$: $\langle \mathscr{S}_{i_n=\alpha}, \mathscr{S}_{i_n=\beta}\rangle = 0$;
\item the ordering: $\norm{\mathscr{S}_{i_n=1}}_F \geq \norm{\mathscr{S}_{i_n=2}}_F \geq \cdots \geq
\norm{\mathscr{S}_{i_n=I_n}}_F \geq 0$ for all $n$.
\end{enumerate}
Given a target multilinear rank $(R_1,\ldots,R_N)$, a truncated HOSVD may be computed in which $\mtx{A}^{(n)}$ contains only $R_n$ orthonormal columns. The procedure is described in Algorithm \ref{alg:HOSVD} \cite{kolda2009tensor}. The important difference between the matrix SVD and the HOSVD is that there is no higher dimensional equivalent of the Eckhart-Young theorem: truncating does not result in an optimal approximation \cite{de2000multilinear}.

\begin{algorithm}[ht]
\caption{Truncated HOSVD}\label{alg:HOSVD}
\begin{algorithmic}[1]
\Procedure{HOSVD}{$\mathscr{X}$,$R_1,\ldots,R_N$}
\For{$n=1,\ldots,N$}
\State $\mtx{A}^{(n)} \gets R_n$ leading left singular vectors of $\mtx{X}_{(n)}$
\EndFor
\State $\mathscr{S} \gets (\mtx{A}^{(1)\intercal},\ldots, \mtx{A}^{(N)\intercal}) \cdot \mathscr{X}$
\State return $\mathscr{S}, \mtx{A}^{(1)},\ldots,\mtx{A}^{(N)}$
\EndProcedure
\end{algorithmic}
\end{algorithm}

The best approximate HOSVD of a given $\mathscr{X}$ is determined by:
\begin{align*}
\min_{\mathscr{S},\{\mtx{A}^{(n)}\}} & \quad \tfrac{1}{2}\norm{\mathscr{X} - (\mtx{A}^{(1)},\ldots, \mtx{A}^{(N)})\cdot\mathscr{S} }_F^2 \\
\text{subject to} & \quad \mathscr{S}\in\mathbb{R}^{R_1\times\cdots\times R_N}, \, \mtx{A}^{(n)}\in\mathbb{R}^{I_n\times R_n}  \text{ and } \mtx{A}^{(n)\intercal}\mtx{A}^{(n)} = \mtx{I}_{R_n},
\end{align*}
which is equivalent to \cite{de2000best}
$$
\begin{aligned}
\max_{\{\mtx{A}^{(n)}\}} & \quad \tfrac{1}{2}\norm{\left(\mtx{A}^{(1)\intercal},\ldots,\mtx{A}^{(N)\intercal}\right)\cdot\mathscr{X}}_F^2 \\
\text{subject to} & \quad \mtx{A}^{(n)}\in\mathbb{R}^{I_n\times R_n} \text{ and } \mtx{A}^{(n)\intercal}\mtx{A}^{(n)} = \mtx{I}_{R_n},
\end{aligned}
$$
where $\mathscr{S} = \left(\mtx{A}^{(1)\intercal},\ldots,\mtx{A}^{(N)\intercal}\right)\cdot\mathscr{X}$. If we use the natural fiber ordering of \cite{elden2009newton} we see that the mode\hyp{}$n$ matricization of $\left(\mtx{A}^{(1)\intercal},\ldots,\mtx{A}^{(N)\intercal}\right)\cdot\mathscr{X}$ is
$$
\mtx{Y}_{(n)} := \mtx{A}^{(n)\intercal}\mtx{X}_{(n)}\left(\mtx{A}^{(1)} \otimes \cdots \otimes \mtx{A}^{(n-1)} \otimes \mtx{A}^{(n+1)} \otimes \cdots \otimes \mtx{A}^{(N)} \right).
$$
Fixing all factor matrices but $\mtx{A}^{(n)}$ gives $\tfrac{1}{2}\norm{\mtx{A}^{(n)\intercal}\mtx{Y}_{(n)}}_F^2$, which is maximized by taking the $R_n$ leading left singular vectors of $\mtx{Y}_{(n)}$ as the columns of $\mtx{A}^{(n)}$. The HOOI Algorithm implements this iteration.

\section{Matrix Manifold Optimization}
\label{appx:manifolds}

In \cite{sterck2016nonlinearly} we discussed how to adapt NPNCG and NPNGMRES to matrix manifolds. In this appendix we summarize the terminology and operations required for optimization on matrix manifolds, using Grassmann and Stiefel manifolds as particular examples, and then describe the NPQN iteration for general matrix manifold problems. 

\subsection{Motivation for Manifold Optimization}
\label{subsect:manifold:motivation}

As previously noted, the tensor Frobenius norm is invariant under orthogonal transformations, meaning \eqref{eq:hosvd:max} does not have isolated maxima. Furthermore, the orthonormality imposed on $\{\mtx{A}^{(n)}\}_{n=1}^N$ introduces a large number of equality constraints. However, if we define the $N$\hyp{}tuple of factor matrices to be a point on a Cartesian product of matrix manifolds, we are able to eliminate both of these issues. This discussion requires some knowledge of two particular types of matrix manifolds, defined below.

The \emph{Stiefel manifold}, $\text{St}(n,p) = \left\{\mtx{X}\in\mathbb{R}^{n\times p} | \mtx{X}^\intercal\mtx{X} = \mtx{I}_p\right\}$, is the set of all $n \times p$ orthonormal matrices. The \emph{Grassmann manifold} (\emph{Grassmannian}), $\text{Gr}(n,p)$, is the set of $p$-dimensional linear subspaces of $\mathbb{R}^n$  \cite{edelman1998geometry}. We can represent $\mathcal{Y}\in\text{Gr}(n,p)$ as the column space of some $\mtx{Y}\in\text{St}(n,p)$. This $\mtx{Y}$ is not unique: the subset of $\text{St}(n,p)$ with the same column space as $\mtx{Y}$ is $\mtx{Y} O_p:= \{\mtx{YM} | \mtx{M}\in O_p\}$, where $O_p$ is the set of $p\times p$ orthogonal matrices. $\text{Gr}(n,p)$ is thus identified with the set of matrix equivalence classes $\text{St}(n,p)/O_p :=\{\mtx{Y}O_p | \mtx{Y}^\intercal\mtx{Y} = \mtx{I}_p\}$ induced by $\mtx{X} \sim \mtx{Y}$ if and only if $\text{Col}(\mtx{X}) = \text{Col}(\mtx{Y})$. The inner product on these manifolds is $\langle \mtx{X}, \mtx{Y} \rangle = \text{tr}(\mtx{X}^\intercal\mtx{Y})$.

If \eqref{eq:hosvd:max} is solved over a Cartesian product of Grassmannians, the representative $N$\hyp{}tuples of factor matrices satisfy the orthonormality constraints by definition. Furthermore, because these matrices represent equivalence classes, the result is an unconstrained problem with isolated extrema:
\begin{equation}
\label{eq:grassopt}
\max_{\{\mtx{A}^{(n)}\}} \quad \tfrac{1}{2}\norm{\left(\mtx{A}^{(1)\intercal},\ldots,\mtx{A}^{(N)\intercal}\right)\cdot\mathscr{X}}_F^2,
\end{equation}
where $\mtx{A}^{(n)}\in\text{St}(I_n,R_n)$ represents $\mathcal{A}^{(n)}\in\text{Gr}(I_n,R_n)$. Expressions for the gradient of this objective function can be found in \cite{elden2009newton, ishteva2011best}.

\subsection{Directions and Movement on Manifolds}
\label{subsect:manifold:movement}

We use \cite{absil2009optimization} as a general reference for this subsection. Let $\mathcal{M}$ denote an arbitrary manifold. A \emph{tangent vector} at $x\in\mathcal{M}$, denoted $\xi_x$, describes a possible direction of travel tangent to $\mathcal{M}$ at $x$. The \emph{tangent space}, $T_x\mathcal{M}$, is the vector space of all tangent vectors at $x$. On $\text{Gr}(n,p)$ a tangent vector is itself an $n\times p$ matrix, and just as $\mtx{Y}\in\text{St}(n,p)$ can represent $\mathcal{Y}\in\text{Gr}(n,p)$, we can use elements of $T_\mtx{Y}\text{St}(n,p)$ to represent elements of $T_\mtx{Y}\text{Gr}(n,p)$. We may express $T_\mtx{Y}\text{St}(n,p)$ as $\mathcal{V}_\mtx{Y} \oplus \mathcal{H}_\mtx{Y}$, where $\mathcal{V}_\mtx{Y}$ contains directions for movement within the equivalence class $\mathcal{Y}$ and $\mathcal{H}_\mtx{Y}$ contains directions for movement between equivalence classes. Elements of $\mathcal{H}_\mtx{Y}$ are used as unique representative tangent vectors for points on $\text{Gr}(n,p)$. Specifically, $\mathcal{H}_\mtx{Y} = \left\{\mtx{Z}\in\mathbb{R}^{n\times p} | \mtx{Y}^\intercal\mtx{Z} = \mtx{0}_p \right\}$, and the orthogonal projection onto $\mathcal{H}_\mtx{Y}$ is \cite{edelman1998geometry, absil2009optimization}
\begin{equation}
\label{eq:orthogonalprojection}
\Pi_\mtx{Y} = \mtx{I} - \mtx{YY}^\intercal.
\end{equation}

Movement along $\mathcal{M}$ in the direction of $\xi$, is described by a \emph{retraction} mapping. The \emph{exponential map} describes motion along a geodesic, the curve connecting two points with minimal length. On $\text{Gr}(n,p)$, the exponential map starting at $\mtx{Y}$ in the direction $\mtx{\xi}$ is
\begin{equation}
\label{eq:exp:retraction}
\text{Exp}_\mtx{Y}(t\mtx{\xi}) = \mtx{YV}\cos(\mtx{\Sigma}t)\mtx{V}^\intercal + \mtx{U}\sin(\mtx{\Sigma}t)\mtx{V}^\intercal,
\end{equation}
where $\xi$ has compact SVD $\mtx{U\Sigma V}^\intercal$. 

Tangent spaces $T_\mathcal{X}\mathcal{M}$ and $T_\mathcal{Y}\mathcal{M}$ for $\mathcal{X}\neq \mathcal{Y}$ are generally different vector spaces, hence linear combinations of $\xi \in T_\mathcal{X}\mathcal{M}$ and $\eta \in T_\mathcal{Y}\mathcal{M}$ are not well-defined. By using a \emph{vector transport} mapping, we instead find a $\xi^\prime \in T_\mathcal{Y}\mathcal{M}$ to use in place of $\xi$. Given $\mathcal{X}\in\text{Gr}(n,p)$ and $\xi,\eta\in T_\mtx{X}\text{Gr}(n,p)$, if $\xi$ has compact SVD $\mtx{U}\pmb{\Sigma}\mtx{V}^\intercal$, the parallel transport of $\eta$ along the geodesic of length $t$ starting at $\mathcal{X}$ in the direction of $\xi$ is \cite{absil2009optimization}:
\begin{equation}
\label{eq:vectortransport}
\mathcal{T}_{\mtx{X},t\xi}(\eta) = \left(
\begin{bmatrix}
\mtx{X}\mtx{V} & \mtx{U}
\end{bmatrix}
\begin{bmatrix}
-\sin(\pmb{\Sigma}t) \\
\cos(\pmb{\Sigma}t)
\end{bmatrix}
\mtx{U}^\intercal
+
(\mtx{I}-\mtx{U}\mtx{U}^\intercal)
 \right)\eta.
\end{equation}
If $\xi = \eta$, this simplifies to
$$
\mathcal{T}_{\mtx{X},t\xi}(\xi) = 
\begin{bmatrix}
\mtx{X}\mtx{V} & \mtx{U}
\end{bmatrix}
\begin{bmatrix}
-\sin(\pmb{\Sigma}t) \\
\cos(\pmb{\Sigma}t)
\end{bmatrix}
\pmb{\Sigma}\mtx{V}^\intercal.
 $$
 If we know that $\mtx{Y} = \text{Exp}_{\mtx{X}}(t\xi)$, we can also compute the approximation $\mathcal{T}_\mtx{Y}(\eta) = \Pi_\mtx{Y}\eta$ using \eqref{eq:orthogonalprojection}.

The direction of travel from $\mathcal{X}$ to $\mathcal{Y}$ cannot be described by a vector $\mathcal{Y}-\mathcal{X}$: this operation is not defined. When $\text{Exp}_\mtx{X}(\cdot)$ is invertible, a tangent vector defining a geodesic from $\mathcal{X}$ to $\mathcal{Y}$ can be found via the \emph{logarithmic map}. Given $\mtx{X}$ and $\mtx{Y}$, the tangent vector in $T_\mtx{X}\text{Gr}(n,p)$ is 
\begin{equation}
\label{eq:logmap}
\text{Log}_\mtx{X}(\mtx{Y}) = \mtx{U}\arctan(\mtx{\Sigma})\mtx{V}^\intercal,
\end{equation}
where $\mtx{U\Sigma V}^\intercal$ is the compact SVD of $\Pi_\mtx{X}\mtx{Y}(\mtx{X}^\intercal\mtx{Y})^{-1}$ \cite{son2013real}.

For a manifold $\mathscr{M} = \prod_{k=1}^{N}\text{Gr}(n_k,p_k)$, a Cartesian product of $N$ Grassmannians, elements are $N$-tuples of linear subspaces $y = (\mathcal{Y}_1,\ldots,\mathcal{Y}_N)^\intercal$, in turn represented by $N$-tuples of matrices $\mtx{y} = (\mtx{Y}_1,\ldots,\mtx{Y}_N)^\intercal$. The tangent space at $y \in \mathscr{M}$ is the Cartesian product of tangent spaces $T_{\mtx{Y}_k}\text{Gr}(n_k,p_k)$. The inner product on $\mathscr{M}$ is $\langle\mtx{x},\mtx{y}\rangle = \sum_{k=1}^N\langle\mtx{X}_k, \mtx{Y}_k\rangle$. All other required operations are performed component\hyp{}wise on $x,y\in\mathcal{M}$, using the operations defined for $\text{Gr}(n_k,p_k)$.

\subsection{NPQN Methods on Matrix Manifolds}
\label{subsect:manifold:NPQN}

To simplify notation bold lowercase letters are used here to represent $n$-tuples of matrices; e.g., $\mtx{x}_k = (\mtx{A}^{(1)}_k, \ldots, \mtx{A}^{(n)}_k)^\intercal$. We define the preconditioner direction to be $\overline{\mtx{g}}_{k} = -\text{Log}_{\mtx{x}_k}(\mathcal{Q}(\mtx{g};\mtx{x}_k))$ (see \eqref{eq:logmap}), the negative of the tangent vector at $\mtx{x}_k$ defining the geodesic to $\mathcal{Q}(\mtx{g};\mtx{x}_k)$. One iteration of manifold NPQN is described in Algorithm \ref{alg:man:NPQN}. The line\hyp{}search is carried out along the curve defined by the retraction $R_{\mtx{x}_k}(\cdot)$ (see \eqref{eq:exp:retraction}), and the vectors $\mtx{s}_k$, $\mtx{y}_k$ and $\overline{\mtx{y}}_k$ require vector transport (see \eqref{eq:vectortransport}). 

\begin{algorithm}[ht]
\caption{NPQN on Manifolds}\label{alg:man:NPQN}
\begin{algorithmic}[1]
\Procedure{manNPQN}{$f, \mtx{g}, \mtx{x}_k$, $\mtx{S}_{k}$, $\mtx{Y}_{k}$, $\overline{\mtx{Y}}_{k}$}
\State $\mtx{g}_k \gets \mtx{g}(\mtx{x}_{k})$
\State $\overline{\mtx{g}}_k \gets -\text{Log}_{\mtx{x}_k}(\mtx{Q}(\mtx{x}_k))$
\State compute $\overline{\mtx{p}}_k$ by \eqref{eq:nonlinear:precond:bfgs:update} or \eqref{eq:nonlinear:precond:broyden:update} \label{line:pk:update}
\State $\mtx{x}_{k+1} \gets R_{\mtx{x}_k}(\alpha_k\overline{\mtx{p}}_k)$ \Comment $\alpha_k$ determined by line\hyp{}search
\State $\mtx{g}_{k+1} = \mtx{g}(\mtx{x}_{k+1})$
\State $\overline{\mtx{g}}_{k+1} \gets -\text{Log}_{\mtx{x}_{k+1}}(\mtx{Q}(\mtx{x}_{k+1}))$
\State $\mtx{s}_k = \mathcal{T}_{\mtx{x}_{k},\alpha_k\overline{\mtx{p}}_k}(\alpha_k\overline{\mtx{p}}_k)$
\State $\mtx{y}_k = \mtx{g}_{k+1} - \mathcal{T}_{\mtx{x}_{k},\alpha_k\overline{\mtx{p}}_k}(\mtx{g}_k)$
\State $\overline{\mtx{y}}_k = \overline{\mtx{g}}_{k+1} - \mathcal{T}_{\mtx{x}_{k},\alpha_k\overline{\mtx{p}}_k}(\overline{\mtx{g}}_{k})$
\State update $\mtx{S}_k$ to $\mtx{S}_{k+1}$, $\mtx{Y}_k$ to $\mtx{Y}_{k+1}$, and $\overline{\mtx{Y}}_k$ to $\overline{\mtx{Y}}_{k+1}$ \label{line:Sk:update}
\State\Return $\mtx{x}_{k+1}, \mtx{S}_{k+1}, \mtx{Y}_{k+1}, \overline{\mtx{Y}}_{k+1}$
\EndProcedure
\end{algorithmic}
\end{algorithm}

In most of this algorithm we work with tangent vectors as $N$\hyp{}tuples of matrices. Two exceptions are in line~\ref{line:pk:update}, where we compute the search direction, and line~\ref{line:Sk:update}, where we update the storage matrices. To use \eqref{eq:nonlinear:precond:bfgs:update} or \eqref{eq:nonlinear:precond:broyden:update} the tangent vectors $\mtx{g}_k$ and $\overline{\mtx{g}}_k$ are converted into 1D arrays by vectorizing each factor matrix column\hyp{}wise and vertically concatenating the results. Once the search direction is computed this process is reversed to produce the tangent vector $\overline{\mtx{p}}_k$ for use in the retraction. Similarly, before the updates in line~\ref{line:Sk:update} the tangent vectors $\mtx{s}_k$, $\mtx{y}_k$, and $\overline{\mtx{y}}_k$ must also be vectorized.

When updating the matrices $\mtx{S}_k$, $\mtx{Y}_k$ and $\overline{\mtx{Y}}_k$, the proper approach for manifold optimization is to parallel transport these matrices to $\mtx{x}_{k+1}$ before appending the new column, as described in \cite[\S~7.1]{savas2010quasi}. However, transport of these matrices can be computationally expensive. Comparisons of algorithms with and without this transport step are recorded in Table \ref{table:tucker:test:m:digit:cvsrch:No:transport}. Based on these results we ultimately decided to omit this step in this paper. 

\end{document}